\documentclass{article}
\usepackage{a4,amsmath,amssymb,amsfonts,amsthm}

\begin{document}
\bibliographystyle{plain}

%
% Macros VM
%
 
% Diff\'erentes fontes
%
\def\M#1{\mathbb#1}     % fontes tableau noir (ensembles R, C, Z, ...)
\def\B#1{\bold#1}       % symboles math\'ematiques accentues
\def\C#1{\mathcal#1}    % lettres caligraphiees (mode math ou texte)
\def\E#1{\scr#1}        % autres fontes caligraphiees (mode math ou texte)
\def\R{\M{R}}           %
\def\Z{\M{Z}}           %
\def\N{\M{N}}           % parce que vraiment on les utilise trop souvent
\def\Q{\M{Q}}           %
%
% Quelques macros d'interet general
%
\def\ds{\displaystyle}
\def\hooklongrightarrow{\lhook\joinrel\longrightarrow}

\def\vfi{\varphi}
\def\epsi{\varepsilon}
\def\op#1{\operatorname{#1}}
\def\SP{\op{Spec}\M{Z}}
\def\ov#1{\overline{#1}}
\def\un#1{\underline{#1}}
\def\Res{\op{Res}}
\def\coef{\op{coef}}
\def\div{\op{div}}
\def\Inv{\op{Inv}}
\def\Spec{\op{Spec}}
\def\im{\op{Im}}
\def\rg{\op{rg}}
\def\Hom{\op{Hom}}
\def\Supp{\op{Supp}}
\def\Arg{\op{Arg}}
\def\bi{\binom}
 \bibliographystyle{plain}
 \def\Aut{{\rm Aut}}
 \def\beginProof{\par{\bf Proof: }}
 \def\endProof{${\qed}$\par\smallskip}
 \def\ar#1{\widehat{#1}}
 \def\mn{{\mu_{n}}}
 \def\pr{^{\prime}}
 \def\prpr{^{\prime\prime}}
 \def\mtr#1{\overline{#1}}
 \def\ra{\rightarrow}
 \def\Bbb{\bf}
 \def\mQ{{\Bbb Q}}
 \def\mR{{\Bbb R}}
 \def\mZ{{\Bbb Z}}
 \def\mC{{\Bbb C}}
 \def\mN{{\Bbb N}}
 \def\Qmn{{\mQ(\mn)}}
 \def\refeq#1{(\ref{#1})}
 \def\umn{^{\mn}}
 \def\lmn{_{\mn}}
 \def\blb{{\big(}}
 \def\brb{{\big)}}
 \def\Hom{{\rm Hom}}
 \def\Tr{{\rm Tr\,}}
 \def\Gal{{\rm Gal}}
 \def\End{{\rm End}}
 \def\Det{{\op{Det}}}
 \def\Td{{\rm Td}}
 \def\ch{{\rm ch}}
 \def\chg{{\rm ch}_{g}}
 \def\torus{{\cal T}}
 \def\Proj{{\rm Proj}}
 \def\Spec{{\rm Spec}\,}
 \def\Qb{\mtr{\Bbb Q}}
 \def\Zn{{{\Bbb Z}/n}}
 \def\deg{{\rm deg}}
 \def\mod{{\rm mod}}
 \def\ac1{\ar{\rm c}_{1}}
 \def\NIm{{\rm Im}}
 \def\NRe#1{{{\rm Re}(#1)}}
 \def\rk{{\op{rk}}}
 \def\Id{{\rm Id}}
 \def\achmn{\ar{\rm ch}_\mn}
 \def\acheq{\ach\lmn}
 \def\ach{\ar{\ch}}
 \def\Tdeq{\Td_{g}}
 \def\Ker{{\rm Ker}}
 \def\HDL{H_{\rm Dlb}}
 \def\HDR{H_{\rm dR}}
 \def\OK{{{\cal O}_K}}
 \def\S1{{S}^{1}}
 \def\ddc{{\mtr{\partial}\partial\over 2\pi i}}
 \def\dgm{{\ar{c}^{\,1}_{\mu_n}}}
 \def\CS{{\rm CS}}
 \def\LS{{\rm LS}}
 \def\mot{{\cal M}}
 \def\prim{{\rm prim}}
 \def\Err{{\rm Err}}
 \def\HK0{H_{K_0}}
 \def\Req{R_g}
 \def\Ind{{\rm Ind}}
 \newcommand{\kkk}[1]{{\Large\bf#1}}
 \def\Alt{{\rm Alt}}
 \def\Per{{\rm Per}}
\def\HM{{\rm HM}}
\def\ACHQ{\ar{\rm CH}_{\mQ(\mn)}}
 \newcommand{\nlabel}[1]{\label{#1}}

 \newtheorem{theor}{Theorem}[section]
 \newtheorem{prop}[theor]{Proposition}
 \newtheorem{cor}[theor]{Corollary}
 \newtheorem{lemma}[theor]{Lemma}
 \newtheorem{sublem}[theor]{sublemma}
 \newtheorem{defin}[theor]{Definition}
 \newtheorem{conj}[theor]{Conjecture}
 
 \newtheorem*{thm1}{Theorem 1}
 \newtheorem*{thm2}{Theorem 2}
 \newtheorem*{coro}{Corollary}

 \author{Vincent MAILLOT\protect\footnote{Institut de Math\'ematiques de Jussieu,
 Universit\'e Paris 7 Denis Diderot,
 Case Postale 7012,
 2 place Jussieu,
 F-75251 Paris Cedex 05, France,
 E-mail : vmaillot@math.jussieu.fr}\ \ and
 Damian ROESSLER\protect\footnote{Mathematik Departement,
 ETH-Zentrum,
 CH-8092 Z\"urich,
 Switzerland}}
 \title{On the periods of motives with complex multiplication and
 a conjecture of Gross-Deligne}
 \maketitle


 \begin{abstract}
 We prove that the existence of an automorphism
 of finite order on a $\Qb$-variety $X$ implies
 the existence of algebraic linear relations between the
 logarithm of certain periods of $X$ and the
 logarithm of special values of the $\Gamma$-function.
 This implies that a slight variation of
 results by Anderson, Colmez and Gross on the periods
 of CM abelian varieties is valid for a larger class of CM motives.
 In particular, we prove a weak form of the period conjecture of
 Gross-Deligne \protect\cite[p. 205]{Gross}\protect\footnote{This 
should not be confused with the conjecture by
 Deligne relating periods and values of $L$-functions.}.
 Our proof relies on the arithmetic fixed point formula
 (equivariant arithmetic Riemann-Roch theorem)
 proved by K. K\"ohler and the second
 author in \protect\cite{LRRI} and the vanishing of the
equivariant analytic torsion for the Dolbeault complex.
 \end{abstract}
 \date
 \parindent=0pt

 \section{Introduction}
 \label{introduction}
 Let $\mot$ be a (homological Grothendieck) motive
 defined over $Q_0$, where $Q_0$ is an algebraic
 extension of $\mQ$ embedded in $\mC$. We shall use the properties
 of the category of motives over a field which 
are listed at the beginning of \cite{Deligne}.  
 The complex singular cohomology 
 $H(\mot,\mC)$ of the manifold of complex points
 of $\mot$ is then endowed with two natural $Q_0$-structures.
The first one is induced by the standard Betti $\mQ$-structure 
$H(\mot,\mQ)$ via the identifications $H(\mot,Q_0) = H(\mot,\mQ)\otimes_{\mQ}Q_0$
and $H(\mot,\mC) = H(\mot,Q_0)\otimes_{Q_0}\mC$
and will be referred to as the Betti (or singular) $Q_0$-structure
on $H(\mot,\mC)$.
The second one arises from the comparison isomorphism between $H(\mot,\mC)$ and
the de Rham-Dolbeault cohomology of $\mot$ (tensorized with $\mC$ over $Q_0$)
and will be referred to as the de Rham $Q_0$-structure.
If $\mot$ is 
a smooth variety $f:Y\ra \op{Spec}Q_0$, its de Rham-Dolbeault cohomology is 
by definition the $Q_0$-vector
space $\HDL(Y):=\oplus_{k\geqslant 0}\oplus_{p+q=k}R^q f_*(\Lambda^p(\Omega(f)))$ which 
is endowed with a natural $\mN$ as well as $\mN\times\mN$ grading. 

Let $Q$ be a finite (algebraic) extension of $\mQ$
and suppose that all the embeddings
of $Q$ into $\mC$ factor through $Q_0$. Furthermore let us suppose
 that $\mot$ is endowed with a $Q$-motive structure (over
 $Q_0$).  A $Q$-motive is also called a motive with
coefficients in $Q$ (see \cite[Par. 2]{Deligne}).
 The $Q$-motive structure
 of $\mot$ induces an inner direct sum of complex vector spaces
 $$
 H(\mot,\mC)=\bigoplus_{\sigma\in\Hom(Q,\mC)}H(\mot,\mC)_\sigma
 $$
 which respects both $Q_0$-structures.
 The notation $H(\mot,\mC)_\sigma$ refers to the complex vector
 subspace of $H(\mot,\mC)$ where $Q$ acts via $\sigma\in\Hom(Q,\mC)$.
 The determinant
 $\det_\mC(H(\mot,\mC)_\sigma)$ thus has
 two $Q_0$-structures. Let $v_{\op{sing}}$ (resp. $v_{\op{dR}}$) be a non-vanishing element
 of $\det_\mC(H(\mot,\mC)_\sigma)$ defined
 over $Q_0$ for the singular (resp. for the de Rham) $Q_0$-structure.
We write $P_\sigma(\mot)$ for the (uniquely defined and independant 
of the choices made) image in 
$\mC^\times/Q_0^\times$ of the complex number
 $\lambda$ such that $v_{\op{dR}}=\lambda\cdot v_{\op{sing}}$.

 Let $\chi$ be an odd simple Artin 
character of $Q$ and let us suppose at this point that $\mot$ 
is homogenous of weight $k$. 
Consider the following conjecture:

 \medskip
 {\bf Conjecture A}$\mathbf{(\boldsymbol{\mot},\boldsymbol{\chi})}$. {\it
 The equality of complex numbers
 \begin{equation*}
 \begin{split}
 \sum_{\sigma\in\Hom(Q,\mC)}
 &\log|P_\sigma(\mot)|^2
 \chi(\sigma) \\
 &=-{L\pr(\chi,0)\over L(\chi,0)}
 \sum_{\sigma\in\Hom(Q,\mC)}
 \sum_{p + q = k}p\cdot\rk
 ({H}^{p,q}(\mot,\mC)_{\sigma})\chi(\sigma)
 \end{split}
 \end{equation*}
 is verified, up to addition of a term of the
form $\sum_{\sigma\in\Hom(Q,\mC)}\log|\alpha_\sigma|\chi(\sigma)$, where
$\alpha_\sigma\in Q_0^\times$.}
\medskip

This conjecture is a slight strengthening of the case $n=1$, 
$Y=\op{Spec}Q_0$ of the conjecture \cite[Conj. 3.1]{MRI}. 
Notice that this conjecture has both a ``motivic'' and an
 ``arithmetic'' content. More precisely, {\it if the Hodge conjecture 
holds} and $Q_0=\Qb$, this conjecture can be reduced to the case where $\mot$ is 
a submotive of an abelian variety with complex multiplication by $Q$.
Indeed, assuming the Hodge conjecture, one can show by examining its associated
 Hodge structures that some exterior
 power of $\mot$ (taken over $Q$) is isomorphic to
 a motive over $\Qb$ lying in the
tannakian category generated by abelian varieties
 with maximal complex multiplication by $Q$.
 In this latter case, the conjecture A
 is contained in a conjecture of Colmez \cite{C}.
Performing this reduction to CM abelian varieties or circumventing 
it is the ``motivic'' aspect of the conjecture. 

However, even in 
the case of CM abelian varieties, the conjecture seems far from proof: 
in the knowledge of the authors, only 
the case of Dirichlet characters has been tackled with up to now; 
obtaining a proof of the conjecture A {\it for non-abelian Artin characters}  
(i.e. for abelian varieties with complex multiplication by a field 
whose Galois group over $\mQ$ is non-abelian) is 
the ``arithmetic'' aspect alluded above.

In this text we shall be concerned with both aspects, 
but our original contribution concerns the ``motivic'' aspect, more precisely 
in finding a way to circumvent the Hodge conjecture. 

We now state a weaker form of the conjecture A. 
 Let $\chi$ be a simple odd Artin character of $Q$ as before,
 and $N$ be a subring of $\Qb$.
Let $\mot_0$ be a motive over $Q_0$ (not necessarily homogenous) and suppose
 that $\mot_0$ is endowed with a $Q$-motive structure (over $Q_0$).
Let $\mot_0^k$ $(k\geqslant 0)$ be the motive corresponding to 
the $k$-th cohomology group of $\mot_0$. 

 \medskip
 {\bf Conjecture B}$\mathbf{(\boldsymbol{\mot_0},N,\boldsymbol{\chi})}$. {\it
 The equality of complex numbers
 \begin{eqnarray*}
 \lefteqn{\sum_{k\geqslant 0}(-1)^k\sum_{\sigma\in\Hom(Q,\mC)}
 \log|P_\sigma(\mot_0^k)|^2
 \chi(\sigma)}\\
 &=&-\sum_{k\geqslant 0}(-1)^k{L\pr(\chi,0)\over L(\chi,0)}
 \sum_{\sigma\in\Hom(Q,\mC)}
 \sum_{p+q=k}p\cdot\rk
 ({H}^{p,q}(\mot,\mC)_{\sigma})\chi(\sigma)\\
 \end{eqnarray*}
 is verified, up to addition of a term of the 
form $\sum_{\sigma\in\Hom(Q,\mC)}\sum_i
(b_{i,\sigma}\log|\alpha_{i,\sigma}|)\chi(\sigma)$, 
where 
$\alpha_{i,\sigma}\in Q_0^\times$, $b_{i,\sigma}\in N$ and $i$ runs over a 
finite set of indices.}
 \medskip

Note that the conjecture A (resp. B) only depends 
on the vector space $H(\mot,\mC)$ (resp. 
$H(\mot_0,\mC)$), together with its 
Hodge structure (over $\mQ$), its de Rham
$Q_0$-structure and its additional $Q$-structure.
If $V$ is a $\mQ$-vector 
space together with the just described structures on $V\otimes_{\mQ}\mC$
(all of them satisfying the obvious compatibility relations), 
we shall accordingly write 
${\rm A}(V,\chi)$ (resp. ${\rm B}(V,N,\chi)$) for the corresponding 
statement, even if $V$ possibly does not arise from a motive.

 In this article we shall prove the conjecture B 
(and to a lesser extent, part of the conjecture A) for a
 large class of motives, which includes abelian varieties with
 complex multiplication by an abelian extension of $\mQ$,
 without assuming the Hodge conjecture (or any
 other conjecture about motives) to reduce the problem to 
 abelian varieties.
Let us add that even in 
the case of abelian varieties, our method of proof is completly 
different from the existing ones.

 A consequence of our results is that on any $\Qb$-variety $X$,
 the existence of a finite group action implies
 the existence of
 non-trivial algebraic linear relations between the
 logarithm of the periods of the eigendifferentials of $X$ (for the action
 of the group) and the
 logarithm of special values of the $\Gamma$-function (recall that 
they are related to the logarithmic derivatives of Dirichlet $L$-functions 
at $0$ via the Hurwitz formula). 
 More precisely, our results are the following:
\smallskip
 
Let $X$ be a smooth and projective variety
together with an automorphism  $g:X\ra X$ of order $n$, everything being defined 
over a number field $Q_0$. 
Let us denote by $\mu_{n}(\mC)$ (resp. $\mu_{n}(\mC)^{\times}$) the
group of $n$-th roots of unity (resp. the set of primitive $n$-th roots of unity) in $\mC$.
Suppose that $Q_0$ is chosen 
large enough so that it contains $\mQ(\mn)$; 
and let $P_n(T) \in \mQ[T]$ be the polynomial
\[
P_{n}(T) = \sum_{\zeta \in \mu_{n}(\mC)^{\times}}
\prod_{\xi \in \mu_{n}(\mC)\backslash\{\zeta\}}\frac{T - \xi}{\zeta - \xi}.
\]
The submotive $\C{X}(g)=\C{X}(X,g)$ cut out in $X$ by
the projector $P_n(g)$ is endowed by construction with a natural $Q := \mQ(\mn)$-motive structure.
\begin{thm1}
 For all the odd primitive Dirichlet
 characters $\chi$ of $\mQ(\mn)$, 
the conjecture ${\rm B}(\C{X}(g),\mQ(\mn),\chi)$ holds.
% \nlabel{TH1}
 \end{thm1}

Let now $Q$ be a finite abelian extension of $\mQ$ with 
conductor $f_Q$ and
 let $\mot_0$ be the motive associated to an abelian variety
 defined over
 $Q_0$ with (not necessarily maximal) complex multiplication by
 ${\cal O}_Q$. We suppose that the action of ${\cal O}_Q$ is defined over $Q_0$ 
and that $\mQ(\mu_{f_Q})\subseteq Q_0$.  
\begin{thm2}
For all the odd Dirichlet characters $\chi$ of $Q$, 
the assertion ${\rm B}(\mot_0^1,\mQ(\mu_{f_Q}),\chi)$ holds.
% \nlabel{TH2}
 \end{thm2}
As a consequence of the existence of
the Picard variety and of theorems 1 and 2, we get:
\begin{coro}
Let the hypothesies of Theorem 1 hold and suppose also that $X$ is a 
surface. 
For all the odd primitive Dirichlet
 characters $\chi$ of $\mQ(\mn)$, 
the conjecture ${\rm B}(H^{2}(\C{X}(X,g)),\mQ(\mn),\chi)$ holds.
\end{coro}

 Our method of proof relies heavily on
 the arithmetic fixed point formula (equivariant arithmetic
 Riemann-Roch theorem) proved by
 K. K\"ohler and
 the second author in \cite{LRRI}. More precisely,
 we write down the fixed point formula as applied to the Dolbeault complex
 of a variety equipped with the action of a finite group. This
 yields a formula for some linear combinations of
 logarithms of periods of the variety in terms of derivatives
 of (partial) Lerch $\zeta$-functions. Using the Hurwitz formula
 and some combinatorics, we can translate this into
 Theorems 1 and 2. In general the fixed point formula of 
 \cite{LRRI},
 like the arithmetic Riemann-Roch theorem, contains an anomaly term,
 given by the equivariant Ray-Singer analytic torsion, which has proved to be
 difficult to compute explicitly. In the case of the Dolbeault complex, this
 anomaly term vanishes for simple symmetry reasons. It is this
 fact that permits us to conclude.

When $Q_0=\Qb$, $Q$ is an abelian 
extension of $\mQ$ and $\mot_0$ is an abelian variety with 
maximal complex multiplication by $Q$, the assertion ${\rm A}(\mot_0^1,\chi)$ was 
proved by Anderson in \cite{Anderson}, whereas 
 the statement ${\rm A}(\mot_0^1,\chi)$ 
 had already been proven by Gross 
\cite[Th. 3, Par. 3, p. 204]{Gross} in the case where
 $Q$ is an imaginary quadratic extension of $\mQ$, $Q_0=\Qb$ and 
$\mot_0$ is an abelian variety with (not necessarily maximal) 
complex multiplication by $Q$.
 One could probably derive Theorem 2
 from the results of Anderson, 
using the result of Deligne on absolute Hodge 
cycles on abelian varieties \cite{DMOS}   
(proved after the theorem of Gross and inspired by it), 
 which can be used as a substitute of the Hodge conjecture 
in this context. 
In the case where $\mot_0$ is an abelian variety with maximal
 complex multiplication and $Q$ is an abelian
extension of $\mQ$, Colmez \cite{C} proves a 
much more precise version of
 ${\rm A}(\mot_0^1,\chi)$. He uses the N\'eron model of
 the abelian variety to normalise the periods so as to eliminate
 all the indeterminacy and proves an equation similar to
 Theorem 2 for those periods. A slightly weaker form of
 his result (but still much more precise than Theorem 2) can
 also be obtained from the arithmetic fixed point formula, when
 applied to the N\'eron models. This is carried out in \cite{LRRIV}. 
Finally, let us mention that 
when $\mot_0$ is the motive of a CM elliptic curve,
Theorem 2 is just a weak form of the Chowla-Selberg
 formula \cite{CS}. For a historical introduction 
 to those results, see \cite[p. 123--125]{Sch}.

 In the last section of the paper, we compare the conjecture
 A with the period conjecture of Gross-Deligne \cite[Sec. 4, p. 205]{Gross}. This
 conjecture is a translation into the language
 of Hodge structures of a special case of conjecture A,
 with $Q$ an abelian extension of $\mQ$. For example,
 we show the following: Theorem 1 implies that
 if $S$ is a surface defined over $\Qb$ and if $S$ is endowed with an
 action of an automorphism $g$ of finite prime order $p$, then the
natural embedding of
 the Hodge structure $\det_{\mQ(\mu_p)}(H^2(\C{X}(S,g),\mQ))$ into
 $H(\times_{r=1}^{\dim_{\mQ(\mu_p)}H^2(\C{X}(S,g),\mQ)}S,\mQ)$ satisfies a weak form
of the period conjecture. Notice
 that all the cases of the period conjecture which have been verified
 in the litterature correspond to motives cut out in varieties with a prescribed
 geometric structure giving rise to (coarse) moduli spaces 
(abelian varieties, Fermat hypersurfaces, Kuga-Sato
 varieties \dots) whereas this result does not specify a geometric
 structure on $S$. 

 In the light of the application of the arithmetic fixed point formula to
 the conjectures A and B, it would be interesting to investigate whether
 this formula is related to the construction of the cycles
 whose existence (postulated by the Hodge conjecture) would 
 be necessary to reduce the conjecture A to abelian varieties.
\medskip

 {\bf Acknowledgments.} It is a pleasure to thank Y. Andr\'e, P. Colmez,
 P. Deligne and C. Soul\'e for suggestions and interesting
 discussions.

\section{Preliminaries}
\label{preliminaries}

\subsection{Invariance properties of the conjectures}

Let $Q_0$ and $Q$ be number fields taken as in the Introduction, and
let $H$ be a (homogenous) Hodge structure (over $\mQ$). 
The $\mC$-vector space $H_\mC:=H_{\mQ}\otimes_\mQ\mC$ comes with 
a natural $Q_0$-structure given by $H_{\mQ}\otimes_\mQ Q_0$.
Suppose that $H_\mC$ is endowed with another $Q_0$-structure.
The first of these two $Q_0$-structures will be referred to
as the Betti (or singular) one, and the second as 
the de Rham $Q_0$-structure on $H_\mC$.
Suppose furthermore that $H_\mC$ is endowed with an additional $Q$-vector space structure compatible 
with both the Hodge structure and the (Betti and de Rham) $Q_0$-structures.
This $Q$-structure induces an inner direct sum of 
$\mC$-vector spaces $H_{\mC} := \oplus_{\sigma \in \op{Hom}(Q,\mC)}H_{\sigma}$.
Let $V:=\oplus_{\sigma\in\Hom(Q,\mC)}\det_\mC(H_\sigma)$ and 
let $m:=\dim_Q(H)$.  There 
is an embedding $\iota:V\hookrightarrow \otimes_{k=1}^m H_\mC$ 
given by $\iota(\oplus_\sigma v_{1}^\sigma
\wedge\dots\wedge v_m^\sigma):=
\sum_\sigma\Alt(v_1^\sigma\otimes\dots\otimes v_m^\sigma)$. Recall that $\Alt$ 
is the alternation map, described by the formula 
$\Alt(x_1\otimes\dots\otimes x_m):={1\over m!}\sum_{\pi\in{\frak 
S}_m}{\rm sign}(\pi)\pi(x_1\otimes\dots
\otimes x_m)$; here ${\frak S}_m$ is the permutation group 
on $m$ elements and $\pi$ acts on $\otimes_{k=1}^m H_\mC$ by permutation 
of the factors. 
\begin{lemma}
The space $V$ inherits the Hodge structure as well as 
the Betti and de Rham
$Q_0$-structures of $\otimes_{k=1}^m H_\mC$ 
via the map $\iota$.
\label{EmbHodge}
\end{lemma}
\beginProof
The bigrading of $H_\mC$ is described by the weight of $H$ and 
by an action $\upsilon:\mC^\times\ra\End_\mC(H_\mC)$ 
of the complex torus $\mC^\times$, which commutes with complex 
conjugation. The bigrading of $\otimes_{k=1}^m H_\mC$ 
is described by the weight $m\cdot{\rm weight}(H)$ and the 
tensor product action $\upsilon^{\otimes m}:\mC^\times\ra\End_\mC
(\otimes_{k=1}^m H_\mC)$. On the other 
hand we can describe a bigrading on 
each $\det_\mC(H_\sigma)$ by the weight $m\cdot{\rm weight}(H)$ 
and by the exterior product action. The map $\iota$ commutes 
with both actions by construction.\\
To prove that $V$ inherits the Hodge $\mQ$-structure, consider 
that there is an action by $\mQ$-vector space 
automorphisms of $\Aut(\mC)$ on 
$\otimes_{k=1}^m H_\mC$ given by $a((h_1\otimes z_1)\otimes\dots 
\otimes(h_m\otimes z_m)):=
(h_1\otimes a(z_1))\otimes\dots 
\otimes(h_m\otimes a(z_m))$. An element $t$ of 
$\otimes_{k=1}^m H_\mC$ is defined over $\mQ$ (for the Hodge 
$\mQ$-structure) if and only if $a(t)=t$ for all $a\in\Aut(\mC)$. 
For each $\sigma\in\Hom(Q,\mC)$, 
let $b_{1}^\sigma,\dots,b_m^\sigma$ be a 
a basis of $H_\sigma$, which is defined over $\sigma(Q)$ and such 
that $a(b_i^\sigma)=b_i^{a(\sigma)}$ for all $a\in\Aut(\mC)$. 
This can be achieved by taking the conjugates under the action 
of $\Aut(\mC)$ of a given basis.   
Now choose a basis $c_{1},\dots,c_{d_Q}$ of $Q$ over $\mQ$ and let
$e_i:=\sum_\sigma\sigma(c_i)b_{1}^\sigma\wedge\dots\wedge b_m^\sigma$. 
By construction, the elements $\iota(e_1),\dots,\iota(e_{d_Q})$ 
are invariant under $\Aut(\mC)$ and they are linearly 
independent over $\mC$, because the determinant 
of the transformation matrix from the 
basis $\{b_1^\sigma\wedge\dots\wedge b_m^\sigma\}_\sigma$ 
to the basis formed by the $e_i$ is the discriminant of 
the basis $e_i$ over $\mQ$. They 
thus define over $V$ a $\mQ$-structure $V_{\mQ}$ which is compatible 
with the Hodge $\mQ$-structure of $\otimes_{k=1}^m H_\mC$.
The Betti $Q_{0}$-structure on $V$ is then just taken to be 
$V_{\mQ}\otimes_{\mQ}Q_{0}$.\\
To show that $V$ inherits the de Rham $Q_0$-structure of $\otimes_{k=1}^m 
H_\mC$, 
just notice that for each $\sigma\in\Hom(Q,\mC)$, the space $H_\sigma$
as a basis $\alpha_1^\sigma,\dots,\alpha_m^\sigma$ 
defined over the de Rham $Q_0$-structure of 
$H_\mC$. The elements $\alpha_1^\sigma\wedge\dots\wedge\alpha_m^\sigma$ 
form a basis of $V$ and $\iota(\alpha_1^\sigma\wedge\dots\wedge\alpha_m^\sigma)$ 
is by construction defined over $Q_0$.
\endProof
In view of the last lemma the complex vector space 
$V$ arises from a (homogenous) Hodge structure over $\mQ$ that we shall 
denote by $\det_Q(H)$.
The embedding $\iota$ arises from 
an embedding of Hodge structures $\det_Q(H)\hookrightarrow\otimes_{k=1}^m H$ 
and $\det_Q(H)$ inherits a Betti and a de Rham $Q_0$-structure from this embedding. 
If $H_{0} = \oplus_{w \in \mZ}H_{w}$ is a direct sum 
of homogenous Hodge structures (graded by the weight), each of them satisfying the hypotheses of lemma \ref{EmbHodge},  
we extend the previous definition to $H_{0}$ by letting $\det_Q(H_{0}) := 
\oplus_{w \in \mZ}\det_Q(H_{w})$.
\begin{prop}
 The assertion ${\rm A}(\mot,\chi)$ (resp. ${\rm B}(\mot_0,N,\chi)$) is equivalent
 to the assertion ${\rm A}(\det_Q(H(\mot,\mQ)),\chi)$ (resp. ${\rm B}(\det_Q(H(\mot_0,\mQ)),N,\chi)$).
\label{DetInv}
 \end{prop}
 \beginProof
 We examine both sides of the equality in
 the assertion ${\rm A}(H(\mot,\mQ),\chi)$, when 
$H(\mot,\mQ)$ is
 replaced by $\det_Q(H(\mot,\mQ))$. 
 From the definition of $\det_Q(H(\mot,\mQ))$, we see that 
 the lefthand side is unchanged. As to the righthand side,
 it is sufficient to show that
 $$
 \sum_{p+q=k}p\cdot\rk(H^{p,q}_\sigma)=
 \sum_{p+q=r\cdot k}p\cdot\rk({\op{det}_{\mC}}{(H_\sigma)^{p,q}}),
 $$
 where $r:=\rk(H_\sigma)$, $k$ is the weight of $\mot$ and $H:=H(\mot,\mQ)$.
 To prove it, we let
 $v_1,\dots, v_r$ be a basis of $H_\sigma$, which
 is homogenous for the grading. The last equality can
 then be rewritten as
 $$
 \sum_{j=1}^r p_{\cal H}(v_j)=p_{\cal H}(v_1\wedge\dots\wedge v_r)
 $$
 (where $p_{\cal H}$ stands for the Hodge $p$-type)
 which holds from the definitions.
 The proof of the second equivalence runs along the 
 same lines.
 \endProof
Let $\mot$, $\mot\pr$ be two $Q$-motives (over $Q_0$). We 
recall that the {\it tensor product with coefficients in $Q$} 
 of $\mot$ and $\mot\pr$, written 
$\mot\otimes_Q\mot\pr$, is the 
the largest direct factor of the motive $\mot\otimes\mot\pr$ 
on which the left and the right action of $Q$ on 
$\mot\otimes\mot\pr$ coincides (see \cite[Par. 2, p. 320]{Deligne}).
If $E$ is a finite extension of $Q$ and $\chi$ is a character of
$Q$, recall that $\Ind^{E}_{Q}(\chi)$ is the character on $E$ (the
 induced character) defined by the formula
 $\Ind^{E}_{Q}(\chi)(\sigma_E):=\chi(\sigma_{E}|_Q)$.

\begin{prop}
 Let $E$ be a finite extension of $Q$, such that all
 the embeddings of $E$ are contained in $Q_0$.
 The statement ${\rm A}(\mot\otimes_Q E,\Ind^{E}_{Q}(\chi))$ 
 (resp. ${\rm B}(\mot_{0}\otimes_Q E,N,\Ind^{E}_{Q}(\chi))$) holds
 if and only if ${\rm A}(\mot,\chi)$ (resp. ${\rm B}(\mot_{0},N,\chi)$) holds.
\label{BasInv}
 \end{prop}
 \beginProof
 Let $r$ be the dimension of $E$ over $Q$.
 The choice of a basis $x_{1},\dots,x_r$ of $E$ as a $Q$-vector space
 induces an isomorphism of 
$Q$-motives $\mot\otimes_Q E\simeq\oplus_{j=1}^r\mot$ 
and thus an isomorphism 
 of $\mC$-vector spaces
 $$
 \bigoplus_{j=1}^r H(\mot,\mC)\simeq
 H((\mot\otimes_{Q}E),\mC)
 $$
 which respects the Hodge structure and both
 $Q_0$-structures.
 Under this isomorphism, we also have a decomposition
 $$
 \bigoplus_{j=1}^r H(\mot,\mC)_{\sigma_Q}\simeq
 \bigoplus_{\sigma_E|\sigma_Q}H((\mot\otimes_Q E),\mC)_{\sigma_E}
 $$
 where $\sigma_Q\in\Hom(Q,\mC)$ and the $\sigma_E\in\Hom(E,\mC)$
 restrict to $\sigma_Q$. This decomposition again
 respects the Hodge structure and both $Q_0$-structures. We now compute
 the lefthand side of the equality predicted 
 by ${\rm A}(\mot_E:=\mot\otimes_Q E,\Ind^{E}_{Q}(\chi))$
 \begin{align*}
 \sum_{\sigma_E}\log |P_{\sigma_E}(&\mot{}_E)|^{2}\, \Ind^{E}_{Q}(\chi)(\sigma_E)\\
 &=\sum_{\sigma_Q}\chi(\sigma_{Q})\sum_{\sigma_E|\sigma_Q}
 \log|P_{\sigma_E}(\mot_E)|^2\\
 &=
 \sum_{\sigma_Q}\chi(\sigma_{Q})\sum_{j=1}^r
 \log|P_{\sigma_Q}(\mot)|^2=
 r\cdot\sum_{\sigma_Q}\log|P_{\sigma_Q}(\mot)|^2\,\chi(\sigma_Q)\\
 \intertext{As for the righthand side, we compute}
 \sum_{\sigma_E}\sum_{p,q}p\cdot\rk(H^{p,q}&
 (\mot_E,\mC)_{\sigma_E})\Ind^{E}_{Q}(\chi)(\sigma_E)\\
 &=
 \sum_{\sigma_Q}
 \sum_{p,q}p\cdot\chi(\sigma_{Q})\,\rk(\oplus_{\sigma_E|\sigma_Q}
 H^{p,q}
 (\mot_E,\mC)_{\sigma_E})\\
 &=
 \sum_{\sigma_Q}
 \sum_{p,q}p\cdot\chi(\sigma_{Q})\cdot r\cdot\rk(
 H^{p,q}
 (\mot,\mC)_{\sigma_Q})
\end{align*}

 dividing both sides by $r$, we are reduced to the conjecture
 ${\rm A}(\mot,\chi)$. The proof of the second equivalence
 is similar.
 \endProof
 
\subsection{The arithmetic fixed point formula}

For the sake of completness and in order to fix notations,
we shall review in this section the arithmetic fixed point formula proved
 by K. K\"ohler and the second author in \cite{LRRI}.
 Many results will be stated without proof;
we refer to \cite[Sec. 4]{LRRI} for more details and 
further references to the litterature.

 Let $D$ be a regular arithmetic ring, i.e. a regular,
 excellent, Noetherian integral ring, together
 with a finite set $\cal S$ of injective ring homomorphisms of $D\hookrightarrow{\Bbb 
C}$,
 which is invariant under complex conjugation.
% \cite[Def. 3.1.1, p. 124]{GS2}
 Let $\mn$ be the diagonalisable group scheme over $D$
 associated to the group $\Zn$.
 An {\bf equivariant arithmetic variety} $f:Y\ra \op{Spec} D$
 is a regular integral scheme,
 endowed with a $\mn$-action over $\op{Spec} D$, such that there exists
 a $\mn$-equivariant ample line bundle on $Y$.
 We write
 $Y({\Bbb C})$ for the complex manifold
 $\coprod_{\sigma\in{\cal S}}Y\otimes_{\sigma(D)}{\Bbb C}$.
 The group $\mn(\mC)$ acts
 on $Y({\Bbb C})$ by holomorphic automorphisms and we shall
 write $g$ for the automorphism corresponding to a fixed primitive
 $n$-th root of unity
 $\zeta=\zeta(g)$.
 The subfunctor of fixed points of the functor
 associated to $Y$ is representable and we
 call the representing scheme the {\bf fixed point scheme} and denote
 it by
 $Y_{\mn}$. It is regular and there are
 natural isomorphisms of complex manifolds
 $Y_{\mn}({\Bbb C})\simeq Y({\Bbb C})_{g}$, where
 $Y({\Bbb C})_{g}$ is the set of fixed points of $Y$ under the action of 
$g$.
 We write $f\umn$ for the map $Y\lmn\ra\op{Spec} D$
 induced by $f$.
 Complex conjugation of coefficients induces an antiholomorphic
 automorphism of $Y({\Bbb C})$ and $Y_{\mu_{n}}({\Bbb C})$, both of
 which we denote by $F_{\infty}$.
 We write
 $\widetilde{\frak A}(Y_{\mu_{n}})$ for
 $\widetilde{\frak A}(Y({\Bbb C})_{g}):=\bigoplus_{p\geqslant 0}({\frak
 A}^{p,p}(Y({\Bbb
 C})_{g})/({\rm Im}\,\partial +{\rm Im}\,\mtr{\partial}))$, where ${\frak
 A}^{p,p}(\cdot)$ denotes the set of smooth complex differential forms 
$\omega$
 of type $(p,p)$
 such that $F_{\infty}^{*}\omega=(-1)^{p}\omega$.\\
 A hermitian equivariant sheaf (resp. vector bundle) on $Y$ is
 a coherent sheaf (resp. a vector bundle) $E$ on $Y$, assumed
 locally free on $Y({\Bbb C})$,
 endowed with a $\mu_{n}$-action which lifts the action of
 $\mu_{n}$ on $Y$ and
 a hermitian metric $h$ on $E_{\Bbb C}$, the bundle associated
 to $E$ on the complex points, which is invariant under
 $F_{\infty}$ and $\mn$. We shall write $(E,h)$ or
 $\mtr{E}$ for
 an hermitian equivariant sheaf (resp. vector bundle). There is
 a natural $({\Bbb Z}/n)$-grading
 $E|_{Y_{\mn}}\simeq\oplus_{k\in{\Bbb Z}/n}E_{k}$ on
 the restriction of $E$ to $Y_{\mn}$, whose terms are orthogonal, because
 of the invariance of the metric. We write $\mtr{E}_{k}$ for the $k$-th
 term $(k\in{\Bbb Z}/n)$, endowed with the induced metric.
 We shall also write $\mtr{E}_{\not = 0}$ for
 $\oplus_{k\in({\Bbb Z}/n)\backslash \{0\}}\mtr{E}_{k}$. 
 
 We write $\ch_g(\mtr{E}):=\sum_{k\in\Zn}\zeta(g)^k\ch(\mtr{E}_k)$
 for the equivariant Chern character form
 $\ch_g(E_{\Bbb C},h)$ associated to the
 restriction of $(E_{\Bbb C},h)$ to $Y_{\mn}({\Bbb C})$. Recall
 also that $\Td_{g}(\mtr{E})$ is the differential form
 $\Td(\mtr{E}_0)\Big(\sum_{i\geqslant
 0}(-1)^{i}\ch_{g}(\Lambda^{i}(\mtr{E}_{\not = 0}))\Big)^{-1}$.
 If ${\cal E}:0\ra E\pr\ra E\ra E\prpr\ra 0$ is an exact sequence
 of equivariant sheaves (resp. vector bundles), we shall write $\mtr{\cal 
E}$
 for the
 sequence
 $\cal E$ together with $\mn(\mC)$ and $F_{\infty}$-invariant hermitian
 metrics on
 $E\pr_{\Bbb C}$,
 $E_{\Bbb C}$ and
 $E\prpr_{\Bbb C}$. To $\mtr{\cal E}$ and $\ch_g$ is associated an 
equivariant
 Bott-Chern secondary class $\widetilde\ch_{g}(\mtr{\cal 
E})\in\widetilde{\frak
 A}(Y_{\mn})$,
 which satisfies
 the equation $\ddc\widetilde\ch_{g}(\mtr{\cal E})=\ch_g(\mtr{E}\pr)+
 \ch_g(\mtr{E}\prpr)-\ch_g(\mtr{E})$.
 \begin{defin}
 The arithmetic equivariant Grothendieck
 group $\ar K^{\mn'}_{0}(Y)$ (resp.
 $\ar{K}^{\mu_{n}}_{0}(Y)$) of $Y$ is the
 free abelian group generated by the elements of
 $\widetilde{\frak A}(Y_{\mu_{n}})$ and by
 the equivariant isometry classes of hermitian equivariant sheaves
 (resp. vector bundles), together with the
 relations
 \begin{description}
 \item[(i)] for every exact sequence $\mtr{\cal E}$ as above,
 $\widetilde\ch_{g}(\mtr{\cal E})=\mtr{E}\pr-\mtr{E}+\mtr{E}\prpr$;
 \item[(ii)] if $\eta\in \widetilde{\frak A}(Y_{\mu_{n}})$ is the sum in
 $\widetilde{\frak A}(Y_{\mu_{n}})$ of
 two elements $\eta\pr$ and $\eta\prpr$, then $\eta=\eta\pr+\eta\prpr$ in
 $\ar K^{\mn'}_{0}(Y)$ (resp. $\ar{K}^{\mu_{n}}_{0}(Y)$).
 \end{description}
 \nlabel{defAEGG}
 \end{defin}
 We shall now define a product on $\ar{K}_{0}^{\mn'}(Y)$ (resp.
 $\ar{K}_{0}\umn(Y)$). Let
 $\mtr{V}$, $\mtr{V}\pr$ be hermitian equivariant sheaves (resp. vector 
bundles)
 and let $\eta,\eta\pr$ be
 elements of $\widetilde{\frak A}(Y_{\mn})$. We define a product $\cdot$ on 
the
 generators of $\ar{K}_{0}^{\mn'}(Y)$
 (resp. $\ar{K}_{0}\umn(Y)$) by the rules
 $
 \mtr{V}\cdot\mtr{V}\pr:=\mtr{V}\otimes\mtr{V}\pr
 $,\
 $
 \mtr{V}\cdot\eta=\eta\cdot\mtr{V}:=\ch_g(\mtr{V})\wedge\eta
 $
 and
 $
 \eta\cdot\eta\pr:=\ddc\eta\wedge\eta\pr
 $ and we extend it by linearity. This product is compatible
 with the relations defining $\ar{K}_{0}^{\mn'}(Y)$
 (resp. $\ar{K}_{0}\umn(Y)$) and defines a commutative ring structure on
 $\ar{K}_{0}^{\mn'}(Y)$ (resp. $\ar{K}_{0}\umn(Y)$).\\
Suppose now that $f$ is proper.  
Fix a $F_{\infty}$-invariant K\"ahler
 metric on $Y({\Bbb C})$, with K\"ahler form $\omega_{Y}$ and suppose
 that $\mn(\mC)$ acts by isometries with respect to this K\"ahler metric.
 Let $\mtr{E}:=(E,h)$ be an equivariant hermitian sheaf on $Y$.
 We write $T_{g}(\mtr{E})$ for the equivariant analytic torsion
 $T_{g}(E_{\Bbb C},h)\in{\Bbb C}$ of
 $(E_{\Bbb C},h)$ over $Y({\Bbb C})$; see \cite[Sec. 2]{K1} or subsection \ref{vanishing} for the definition.
 Let $f:Y\ra \op{Spec} {D}$ be the structure morphism.
 We let $R^{i}f_{*}\mtr{E}$ be the $i$-th direct image sheaf of $E$ 
endowed with its
 natural equivariant structure and $L^{2}$-metric (see \cite[Par. 1.1, p. 22]{GS2} 
for the definition of the latter).
 We also write $\mtr{H^{i}(Y,E)}$ for
 $R^{i}f_{*}\mtr{E}$ and
 $R^{\cdot}f_{*}\mtr{E}$ for the linear combination
 $\sum_{i\geqslant 0}(-1)^{i}R^{i}f_{*}\mtr{E}$.
 Let $\eta\in\widetilde{\frak A}(Y_{\mn})$ and
 consider the rule which associates the element
 $R^{\cdot}f_{*}\mtr{E}-T_{g}(\mtr{E})$ of
 $\ar K^{\mn'}_{0}({D})$ to $\mtr{E}$ and the element
 $\int_{Y({\Bbb C})_{g}}\Td_{g}(\mtr{TY})\eta\in\ar K^{\mn'}_{0}({D})$ to
 $\eta$.
 \begin{prop}
 The above rule descends to a well defined group homomorphism
 $f_{*}:\ar K^{\mn'}_{0}(Y)\ra \ar K^{\mn'}_{0}(D)$.
 \nlabel{WFProp}
 \end{prop}
 One can show that $\ar{K}_0\umn(D)$ is isomorphic to
 $\ar{K}_0^{\mn'}(D)$ via the natural map so that by composition
 the last proposition yields a map
 $\ar{K}^{\mu_{n}}_{0}(Y)\ra \ar{K}^{\mu_{n}}_{0}(D)$, which we
 shall also call $f_{*}$.
 
 Finally, to formulate the fixed point theorem, we define
 the homomorphism
 $\rho:\ar{K}^{\mu_{n}}_{0}(Y)\ra\ar{K}^{\mu_{n}}_{0}(Y_{\mu_{n}})$, which
 is obtained by restricting all the involved objects from
 $Y$ to $Y_{\mu_{n}}$. If $\mtr{E}$ is a hermitian vector bundle
 on $Y$, we write
 $\lambda_{-1}(\mtr{E}):=\sum_{k=0}^{{\rm rk}(E)}(-
1)^{k}\Lambda^{k}(\mtr{E})
 \in\ar{K}^{\mu_{n}}_{0}(Y)$, where $\Lambda^{k}(\mtr{E})$ is the
 $k$-th exterior power of $\mtr{E}$, endowed with its natural
 hermitian and equivariant structure. If $\mtr{E}$ is the orthogonal direct
 sum of two hermitian
 equivariant vector bundles $\mtr{E}\pr$ and $\mtr{E}\prpr$, then
 $\lambda_{-1}(\mtr{E})=\lambda_{-1}(\mtr{E}\pr)\cdot\lambda_{-
1}(\mtr{E}\prpr)$.
 Let $R(\mu_{n})$ be
 the Grothendieck group of finitely generated
 projective $\mn$-comodules over $D$. There are
 natural isomorphisms $R(\mu_{n})\simeq {K_{0}(D)}[{\Bbb Z}/n]
 \simeq{K_{0}(D)}[T]/(1-T^{n})$. Let
 $\mtr{I}$ be the $\mn$-comodule whose
 term of degree $1$ is $D$ endowed with the trivial metric and whose
 other terms are $0$. We make
 $\ar{K}^{\mu_{n}}_{0}(D)$ an $R(\mu_{n})$-algebra under
 the ring morphism which sends $T$ to $\mtr{I}$.
 In the next theorem, which is the arithmetic fixed point
 formula, let
 $\cal R$ be any $R(\mn)$-algebra such that the elements
 $1-T^{k}$ ($k=1,\dots,n-1$) are invertible in $\cal R$.
 The algebra which is minimal with respect to this property
 is the ring
 $R(\mn)_{\{1-T^{k}\}_{k=1,\dots ,n-1}}$, the localization
 of $R(\mn)$ at the multiplicative subset
 generated by the elements $\{1-T^{k}\}_{k=1,\dots ,n-1}$.
 %If $D={\Bbb Z}$, we can make the complex numbers
 %$\bf C$ an $R(\mu_{n})$-algebra under the ring morphism which sends $T$ to an
 %$n$-th root of unity; this gives a possible choice of
 %$\cal R$ in this case.

 Let now
 $\theta\in {\bf R}$. For all $s\in{\Bbb C}$ such that $\Re(s)>1$ we define
Lerch's partial $\zeta$-functions
 $
 \zeta(\theta,s):=\sum_{n\geqslant 1}{\cos(n\theta)\over n^{s}}
 $
 and
 $
 \eta(\theta,s):=\sum_{n\geqslant 1}{\sin(n\theta)\over n^{s}}
 $, and using analytic continuation, we extend them 
to meromorphic functions of $s$ over $\Bbb C$.
 Let $R(\theta,t)$ be the formal power series
 $$
 \sum_{n\geqslant 1, \text{ $n$ {\bf odd}}}(2\zeta\pr(\theta,-n)+
 \sum_{j=1}^{n}{\zeta(\theta,-n)\over j}){t^{n}\over n!}+
 i\sum_{n\geqslant 0, \text{ $n$ {\bf even}}}(2\eta\pr(\theta,-n)+
 \sum_{j=1}^{n}{\eta(\theta,-n)\over j}){t^{n}\over n!}.
 $$
 We shall need $R(\theta,(\cdot))$ which
 is the unique additive characteristic class on holomorphic
 vector bundles such that
 $R(\theta,L)=R(\theta,c_{1}(L))$
 for each line bundle $L$.
 Let $V$ be a $\mn$-equivariant vector bundle on $Y$;
 we define
 $$
 R_{g}(V):=\sum_{k=1}^{\op{rk}(V)}R(\arg(\zeta(g)^k),V_{k}).
 $$
 In the next theorem, we consider that the values of $R_g$ lie in
 $\widetilde{\frak A}(Y_{\mn})$.
 \begin{theor}
 Let $\mtr{N}_{Y/Y_{\mu_{n}}}$ be the normal
 bundle of $Y_{\mu_{n}}$ in $Y$, endowed with its quotient equivariant
 structure
 and quotient metric structure (which is $F_{\infty}$-invariant).
 \begin{description}
 \item[(i)] The element $\Lambda:=\lambda_{-
1}(\mtr{N}_{Y/Y_{\mu_{n}}}^{\vee})$
 has an inverse in $\ar{K}^{\mu_{n}}_{0}(Y_{\mu_{n}})\otimes_{R(\mn)}{\cal 
R}$;
 \item[(ii)] let $\Lambda_{R}:=\Lambda\cdot(1+R_{g}(N_{Y/Y\lmn}))$; 
the diagram\footnote{Note that a misprint found its way 
into the statement \cite[Th. 4.4]{LRRI} ($=$ Th. \ref{PMainTh}).
In \cite[Th. 4.4]{LRRI} the term $\Lambda\cdot(1-R_{g}(N_{Y/Y\lmn}))$ has to be replaced 
by $\Lambda\cdot(1+R_{g}(N_{Y/Y\lmn}))$ in the expression for $\Lambda_{R}$.}
 $$
 \begin{matrix}
 \ar{K}^{\mu_{n}}_{0}(Y)&
 \stackrel{\Lambda_{R}^{-1}\cdot\rho}{\longrightarrow}
 &\ar{K}^{\mu_{n}}_{0}(Y_{\mu_{n}})\otimes_{R(\mu_{n})}{\cal R}\cr
 \downarrow\ f_{*} & &\downarrow\ f^{\mu_{n}}_{*}\cr
 \ar{K}^{\mu_{n}}_{0}({D})&
 \stackrel{{\rm Id}\otimes 1}{\longrightarrow}
 &\ar{K}^{\mu_{n}}_{0}({D})\otimes_{R(\mu_{n})}{\cal R}\cr
 \end{matrix}
 $$
 commutes.
 \end{description}
 \nlabel{PMainTh}
 \end{theor}
 The proof of this theorem is the object of \cite{LRRI}, it
 combines the deformation to the normal cone technique with
 deep results of Bismut on the behaviour of equivariant analytic
 torsion under immersions \cite{B1}.

\subsection{The equivariant analytic torsion and the $L^2$-metric of 
the Dolbeault complex}
\label{vanishing}

In this subsection, we shall prove the vanishing of the 
equivariant analytic torsion for the Dolbeault complex. 
Before doing so, we shall review some 
results on the polarisation induced by an ample line bundle
on the singular cohomology of a complex manifold.
 
 So let $M$ be a complex projective manifold of dimension $d$ and
 $L$ be an ample line bundle on $M$. Let us
 denote by $\omega\in H^{2}(M,\mC)$ the first Chern class of $L$ and
 for $k\leqslant d$, let us note
 $P^{k}(M,\mC)\subseteq H^{k}(M,\mC)$ the
 primitive cohomology associated to $\omega$; this is a Hodge substructure 
 of $H^{k}(M,\mC)$. Recall that for any $k\geqslant 0$,
 the primitive decomposition
 theorem establishes an isomorphism
 $$
 H^k(M,\mC)\simeq\oplus_{r\geqslant
 \max(k-d,0)}\,\omega^{r}\wedge P^{k-2r}(M,\mC).
 $$
% Let $C:H^k(M,\mC)\ra H^k(M,\mC)$ be the Weil operator, which is defined
% by the rule $C(\eta):=\sum_{p,q\geq 0}i^{p-q}\eta^{p,q}$, where 
% $\eta^{p,q}$ is the component of type $p,q$ of $\eta$ for 
% the Hodge decomposition of $H^k(M,\mC)$. 
 Define the cohomological star operator $*:H^{k}(M,\mC)\ra
 H^{2d-k}(M,\mC)$ by the rule $*\,\omega^{r}\wedge \phi :=
 i^{p-q}(-1)^{(p+q)(p+q+1)/2}{r!\over
 (d-p-q-r)!}\omega^{d-p-q-r}\wedge \phi)$ if $\phi$ is a primitive
 element of Hodge pure type $(p,q)$ and extend it by additivity. We can now
 define a hermitian metric on $H^{k}(M,\mC)$ by the
 formula
 $$
 (\nu,\eta)_L:=\int_{M}\nu\wedge *\,\mtr{\eta}
 $$
 for any $\nu,\eta\in H^{k}(M,\mC)$. This metric is sometimes called the 
 Hodge metric. The next lemma follows from the definition
 of the $L^2$-metric,
% resultat de Weil est aussi dans Wells, th. 3.16, p. 187 
Hodge's theorem on the
 representability of cohomology classes by harmonic forms and the Hodge-K\"ahler 
identities.
 \begin{lemma}
\label{HDRLem}
 Endow $M$ with a K\"ahler metric whose K\"ahler
 form represents the cohomology class of $\omega$ and equip 
the bundles $\Lambda^p(\Omega_M)$ with the corresponding metrics. 
Endow $\oplus_{p+q=k}
 H^{q}(\Lambda^{p}({\Omega}_M))$ with the $L^2$-metric and 
$H^k(M,\mC)$ with the Hodge metric.
 The Hodge-de Rham isomorphism $H^k(M,\mC)\simeq\oplus_{p+q=k}
 H^{q}(\Lambda^{p}({\Omega}_M))$ is an isometry. 
 \end{lemma}
 
 Suppose now that $M$ is a complex compact K\"ahler manifold 
 endowed with a unitary automorphism $g$, and
 let $\mtr{E}$ be a hermitian holomorphic vector bundle on
 $M$ which is equipped with a unitary lifting of the action of $g$.
 Let $\square_{q}^{E}$ be the differential operator
 $(\mtr{\partial}+\mtr{\partial}^*)^{2}$ acting on the
  $C^\infty$-sections of the bundle
 $\Lambda^{q}T^{*(0,1)}M\otimes E$. This space of
 sections is equipped with the $L^{2}$-metric and the
 operator $\square_q^{E}$ is symmetric for that metric; we let
 ${\rm Sp}(\square_q^{E})\subseteq\mR$ be the set of eigenvalues
 of $\square_q^E$ (which is discrete and bounded 
from below) and we let $\op{Eig}_{\,q}^{E}(\lambda)$ be
 the eigenspace associated to an eigenvalue $\lambda$ 
(which is finite-dimensional). Define
 $$
 Z(\mtr{E},g,s):=\sum_{q\geqslant 1}(-1)^{q+1}q\sum_{\lambda\in{\rm
 Sp}(\square_{q}^{E})\backslash\{0\}}
 {\rm Tr}(g^*|_{\op{Eig}_{\,q}^{E}(\lambda)})\lambda^{-s}
 $$
 for $\Re(s)$ sufficiently large.
 The function $Z(\mtr{E},g,s)$ has a meromorphic continuation to the whole plane,
 which is holomorphic around $0$ (see \cite{K1}). By definition, 
 the equivariant analytic torsion of $\mtr{E}$ is given by
 $T_g(\mtr{E}) := Z\pr(\mtr{E},g,0)$.
 
 The non-equivariant analog of the following lemma 
 (whose proof is similar) can be found in \cite{RS}.
 \begin{lemma}
 \label{torsion}
 Let $M$ be a complex compact K\"ahler manifold and let $g$
 be a unitary automorphism of $M$. The identity
 $$
 \sum_{p\geqslant 0}(-1)^{p}T_{g}(\Lambda^p(\mtr{\Omega}_{M}))=0
 $$
 holds.
 \end{lemma}
 \beginProof
  Recall the Hodge decomposition
 (see \cite[Chap. IV, no. 3, Cor. 2]{Weil})
 $$
 A^{p,q}(M)={\cal H}^{p,q}(M)\oplus{\partial}(A^{p-1,q}(M))\oplus
 {\partial}^*(A^{p+1,q}(M))
 $$
where ${\cal H}^{p,q}(M)$ are the harmonic forms for 
the usual Kodaira-Laplace operator $\square_q=(\partial+\partial^*)^2=
(\mtr{\partial}+\mtr{\partial}^*)^2$ and 
$A^{p,q}(M)$ is the space of $C^\infty$-differential forms 
 of type $(p,q)$ on $M$. 
 Let us write $A_{1}^{p,q}(M)$ for ${\partial}(A^{p-1,q}(M))$ and
 $A_{2}^{p,q}(M)$ for ${\partial}^*(A^{p+1,q}(M))$.
 The map ${\partial}|_{A_{2}^{p,q}(M))}$ is an injection and its
 image is $A^{p+1,q}_{1}(M)$. Notice also  
 that the operator $\square_q$ commutes with ${\partial}$ 
and ${\partial}^*$. 
Notice as well 
 that the $C^\infty$-sections of
 $\Lambda^{q}T^{*(0,1)}M\otimes\Lambda^{p}(\Omega_{M})$ correspond to the space
 $A^{p,q}(M)$ and that 
$\square_q^{\Lambda^p(\Omega)}=\square_q|_{A^{p,q}}$. 
For $\lambda\in\mR^{\times}$, we write
 $L_{\lambda}^{p,q}=\Ker(\square_{q}^{\Lambda^{p}(\Omega)}-\lambda)$,
 $L_{\lambda,1}^{p,q}=L_{\lambda}^{p,q}\cap A_{1}^{p,q}$ and
 $L_{\lambda,2}^{p,q}=L_{\lambda}^{p,q}\cap A_{2}^{p,q}$. We compute
 \begin{equation*}
 \sum_{p\geqslant 0}(-1)^{p}{\rm Tr}(g^*|_{L_{\lambda}^{p,q}})=
 \sum_{p\geqslant 0}(-1)^{p}[{\rm Tr}(g^{*}|_{L_{\lambda,1}^{p,q}})+
 \Tr(g^{*}|_{L_{\lambda,2}^{p,q}})]=0
 \end{equation*}
 and from this, we conclude that $\sum_{p\geqslant 0}(-
1)^{p}Z(\Lambda^{p}(\mtr{\Omega}_{M}),g,s)\equiv 0$.
 \endProof

\subsection{An invariant of equivariant arithmetic $K_0$-theory}
\label{invariant}

 From now on, we restrict ourselves
 to the case $D=Q_0$, where $Q_0 \overset{\iota_{0}}{\hooklongrightarrow} \mC$ is a number field 
 embedded in $\mC$, and we fix a primitive $n$-th root of unity
 $\zeta:=e^{2\pi i/n}$. We use this choice 
 to identify the set $\mn(\mC)^\times$ of primitive
 $n$-th roots of unity with the Galois group
 $G:=\Gal(\mQ(\mn)/\mQ) = \Hom(\mQ(\mn),\mC)$. The ring morphism
 $R(\mn)\ra\mQ(\mn)$ which 
sends the generator $T$ on $\zeta$
makes $\mQ(\mn)$ an $R(\mn)$-algebra 
and allows us to take ${\cal R}:=\mQ(\mn)$.
We let 
${\widehat{\op{CH}}}(Q_{0})$ be the arithmetic Chow ring of 
$Q_{0}$ with the set of embeddings ${\cal S} := \{\iota_{0}, \overline{\iota}_{0}\}$, 
in the sense of Gillet-Soul\'e (see \cite{GS1}). There 
is a natural isomorphism ${\widehat{\op{CH}}}(Q_{0})\simeq \mZ\oplus\mR/\log|Q_{0}^\times|$  
and a ring isomorphism $\ar{K}_0(Q_{0})\simeq\ar{\rm CH}(Q_{0})$ 
given by the arithmetic Chern character $\ar{\rm ch}$ (see \cite{GS1}), 
the ring $\ar{K}_0(Q_{0})$ being defined 
similarly to the ring $\ar{K}_0^{\mu_1}(Q_{0})$, with 
${\frak A}^{p,p}(\cdot)$ replaced by the space ${\frak A}^{p,p}_\mR(\cdot)$ 
of real (not complex) 
differential forms  of type $(p,p)$. The ring structure 
on $\mZ\oplus\mR/\log|Q_{0}^\times|$ is given by the formula 
$(r\oplus x)\cdot (r\pr\oplus x\pr):=(r\cdot r\pr,r\cdot x\pr+r\pr\cdot x)$. 
On generators of $\ar{K}_0(Q_{0})$, the arithmetic 
Chern character is defined as follows:  
For $\mtr{V}$  
a hermitian vector bundle on
 $\op{Spec}Q_0$, the arithmetic Chern character $\ar{\ch}(\mtr{V})$ is 
the element $\rk(V)\oplus(-\log||s||^2)$, where 
$s$ is a non-vanishing section of $\op{det}(V)$ and $||\cdot||$ is 
the norm on $\op{det}(V)_\mC$ induced by the metric on $V_\mC$. 
For an element $\eta\in{\frak A}_\mR(\op{Spec}Q_0)\simeq\mR$, the 
arithmetic Chern character $\ar{\ch}(\eta)$ is 
the element $0\oplus{1\over 2}
 \eta$. 

Let now $\aleph_{0}$ be the
 additive subgroup of $\mC$ generated by the elements
 $z\cdot\log|q_0|$ where $q_0\in Q_0^\times$ and
 $z\in\mQ(\mn)$. We define $\ACHQ(Q_{0}):=\mQ(\mn)\oplus\mC/\aleph_{0}$ and 
we define a ring structure on $\ACHQ(Q_{0})$ by the rule 
$(z,x)\cdot (z\pr,x\pr):=(z\cdot z\pr,z\cdot x\pr+z\pr\cdot x)$. Notice that 
there is a natural ring morphism $\psi :{\widehat{\op{CH}}}(Q_{0})\ra\ACHQ(Q_{0})$ and 
that there is a natural $\mQ(\mn)$-module structure on 
$\ACHQ(Q_{0})$. Define a rule which associates elements 
of $\ACHQ(Q_{0})$ to generators of $\ar{K}^{\mu_{n}}_{0}(Q_{0})$ 
as follows. Associate 
the element $\zeta^k\cdot \psi(\ar{\ch}(\mtr{V}))$ 
to a $\mn$-equivariant hermitian vector bundle $\mtr{V}$ of 
pure degree $k$ (for the natural $({\Bbb Z}/n)$-grading)
on $\op{Spec}Q_0$; 
furthermore
 associate the element
 $0\oplus{1\over 2}
 \eta$ to $\eta\in\widetilde{\frak A}(\op{Spec}Q_0)\simeq\mC$.
 \begin{lemma}
 The above rule extends to a morphism of
 $R(\mn)$-modules $\achmn:\ar{K}^{\mu_{n}}_{0}(Q_{0})\ra
 \ACHQ(Q_{0})$.
 \nlabel{UnMor}
 \end{lemma}
\beginProof
Let
 \begin{equation*}
 {\cal V}:\ 0\ra V\pr\ra V\ra V\prpr\ra 0
% \nlabel{SeqV}
 \end{equation*}
 be a an exact sequence of $\mn$-equivariant vector bundles 
($\mn$-comodules) 
 over $Q_{0}$. We endow the members of $\cal V$ with
 (conjugation invariant) hermitian metrics $h\pr$, 
$h$ and $h\prpr$ respectively, such
 that the pieces of the various gradings are orthogonal. 
 The equality 
 \begin{equation*}
 \widetilde{\ch}_\zeta(\mtr{\cal V})=
\sum_{k\in\Zn}\zeta^k\widetilde{\ch}(\mtr{\cal V}_{k})
% \nlabel{Eqchg}
 \end{equation*}
holds (see \cite[Th. 3.4, Par. 3.3]{LRRI}).
 From this and the well-definedness of 
the arithmetic Chern character, the result follows. 
\endProof
We shall write $\dgm$ for the second component 
of $\achmn$, i.e. the component lying in $\mC/\aleph_{0}$. 

 \begin{lemma}
 Let $\mtr{V}$ be a hermitian
 $\mn$-equivariant vector bundle on $\op{Spec} Q_0$. The equation
 %$$
 %\prod_{l\in\Zn}(1-\zeta^l)^{-\rk(V_l)}\,\dgm(\lambda_{-
 %1}(\mtr{V}))=
 %-\sum_{l\in\Zn}{\zeta^l\over 1-\zeta^l}\,\dgm(\mtr{V}_l)
 %$$
\begin{equation}
\prod_{l\in\Zn}(1-\zeta^l)^{-\rk(V_l)}\achmn(\lambda_{-1}(\mtr{V}))=
1 \oplus(-\sum_{l\in\Zn}{\zeta^l\over 1-\zeta^l}\,\dgm(\mtr{V}_l))
\label{GRRR}
\end{equation}
 holds.
\label{2Ten}
\end{lemma}
\beginProof
We shall make use of the canonical isomorphism 
$$
\op{det}(\Lambda^k({W}))\simeq\op{det}({W})^{\otimes{(r-1)!\over 
(r-k)!(k-1)!}}
$$
valid for any vector space $W$ of rank $r$ over a field and any $1 \leqslant k \leqslant r$, 
and constructed as follows:
For any basis $b_1,\dots,b_r$ of $W$, the element 
$$
\bigwedge_{1\leqslant i_1<\dots <i_k\leqslant r}(b_{i_1}\wedge\dots\wedge b_{i_k})
$$ 
of $\op{det}(\Lambda^k({W}))$ is sent on the element
$$
\bigotimes_{j=1}^{(r-1)!\over (r-k)!(k-1)!}(b_1\wedge\dots\wedge b_r)
$$
of $\op{det}({W})^{\otimes{(r-1)!\over (r-k)!(k-1)!}}$. 
This isomorphism is by construction invariant under base change 
to a field extension. Furthermore, 
if one applies the above description to the orthonormal 
basis of a vector space over $\mC$ endowed with a 
hermitian metric, one find that this isomorphism 
is also an isometry for the natural metrics on both sides.
Thus, for any hermitian vector bundle $\mtr{W}$ over $\op{Spec}Q_0$,
the isomorphism of vector bundles
$$
\op{det}(\Lambda^k(\mtr{W}))\simeq\op{det}(\mtr{W})^{\otimes{(r-1)!
\over (r-k)!(k-1)!}}
$$
is an isometry.
Using the definition of the ring structure of 
$\ACHQ(Q_{0})$ and the fact that 
$\achmn(\lambda_{-1}(\mtr{V}\oplus\mtr{V}\pr))=
\achmn(\lambda_{-1}(\mtr{V}))\cdot\achmn(\lambda_{-1}(\mtr{V}\pr))$ for 
any two hermitian equivariant vector bundles over $Q_{0}$, we see that 
as functions of $\mtr{V}$, 
both sides of the equality in \refeq{GRRR} are 
multiplicative for direct sums 
of hermitian vector bundles. We are thus reduced to prove 
the equality 
$$
(1-\zeta^l)^{-\rk(V_l)}\achmn(\lambda_{-1}(\mtr{V}_l))=
1\oplus{-\zeta^l\over 1-\zeta^l}\,\dgm(\mtr{V}_l)
$$
for all $l\in\Zn$. Let $r_l :=\rk(V_l)$ and let 
us write $\ar{c}_{1}(\cdot)$ for $\ach^{[1]}(\cdot)$; we compute 
\begin{eqnarray*}
\achmn(\lambda_{-1}(\mtr{V}_l))&=&\sum_{k=0}^{r_l}
(-1)^k\zeta^{lk}\ach(\Lambda^k(\mtr{V}_l))\\
%&=&
%\sum_{k=0}^{r_l}(-1)^k\zeta^{kl}
%[{r_l!\over k!(r_l-k)!}\oplus{(r_l-1)!\over (r_l-k)!(k-
%1)!}\ar{c}_{1}(\mtr{V}_l)]\\
&=&
\sum_{k=0}^{r_l}(-1)^k\zeta^{lk}{r_l!\over k!(r_l-k)!}
\oplus(\sum_{k=1}^{r_l}(-1)^k\zeta^{lk}{(r_l-1)!\over 
(k-1)!(r_l-k)!})\ar{c}_{1}(\mtr{V}_l).\\
\end{eqnarray*}
Using the binomial formula, we see that the  
last expression can be rewritten as 
$$
(1-\zeta^l)^{r_l}\oplus(-\zeta^l(1-\zeta^l)^{r_l-1})\ar{c}_{1}(\mtr{V}_l)
$$
and the result follows. 
\endProof

\section{Proof of Theorems 1 and 2}

\subsection{Two lemmas}

 For $z$ belonging to the unit circle $\C{S}_1$,
 we define the Lerch's $\zeta$-function $\zeta_L(z,s):=\sum_{k\geqslant 1}{z^k\over k^s}$
 for $s \in \mC$ such that $\Re(s) > 1$, and using analytic continuation, 
 we extend it to a meromorphic function of $s$ over $\mC$.
 
 \begin{lemma}
 \label{CoHoLe}
 Let $Y$ be a scheme smooth over $\mC$ and
 let $E$ be a vector bundle on $Y$ together with an automorphism
 $g:E\ra E$ of finite order (which leaves the zero-section fixed). Let
 $\kappa$ be the class
 $$
 \kappa:=\Td(E_0){\sum_{p\geqslant 0}(-1)^p p\cdot \chg(\Lambda^p(E^{\vee}))\over
 \sum_{p\geqslant 0}(-1)^p\chg(\Lambda^p(E_{\not=0}^{\vee}))}.
 $$
 The equality
\[
\kappa^{[l + \op{rk}(E_{0})]} = - c^{\op{top}}(E_{0})
\sum_{z \in \C{S}_{1}}
\zeta_{\op{L}}(z,-l)\op{ch}^{[l]}(E_{z}^{\vee})
\]
 holds. 
\end{lemma}
\beginProof
According to the splitting principle, we may suppose that 
$E$ is an equivariant direct sum of 
$r := \op{rk}(E)$ line bundles $F_{i}$ on which $g$ acts 
by multiplication with the eigenvalue $\alpha_{i}^{-1} \in \C{S}_{1}$. 
If we set  
$\gamma_{i} := {c}_{1}(F_{i}^{\vee})$ for all $i$, we can 
write 
\begin{equation}
\label{eq_lemme_1}
\sum_{p \geqslant 0} (-1)^p p \cdot \op{ch}_{g}(\wedge^p(E^{\vee}))
= 
\sum_{p \geqslant 0} (-1)^{p} p \sum_{1 \leqslant i_1 < \dots < i_p \leqslant r}
\alpha_{i_1} \dots \alpha_{i_p} 
e^{\gamma_{i_1} + \dots + \gamma_{i_p}}
\end{equation}
If we take the formal derivative (for $t$) of the identity
\[
\prod_{i=1}^{r} ( 1 - \alpha_{i}e^{\gamma_{i}}t) 
= 
\sum_{p \geqslant 0} (-1)^{p}
\left[ \sum_{1 \leqslant i_{1} < \dots < i_{p} \leqslant r} \alpha_{i_{1}}\dots
\alpha_{i_{p}} e^{\gamma_{i_1} + \dots + \gamma_{i_p}} \right]
t^{p}
\]
set $t=1$ and apply  
(\ref{eq_lemme_1}), we obtain
\begin{equation}
\label{eq_lemme_2}
\sum_{p \geqslant 0} (-1)^p p \cdot \op{ch}_{g}(\wedge^p(E^{\vee})) = 
- \prod_{i=1}^{r}(1 - \alpha_{i}e^{\gamma_{i}}) \sum_{j=1}^{r}
\frac{\alpha_{j}e^{\gamma_{j}}}{1 - \alpha_{j}e^{\gamma_{j}}}\, .
\end{equation}
Notice now that we can write
\begin{equation*}
\begin{split}
\prod_{i=1}^{r}(1 - \alpha_{i}e^{\gamma_{i}}) &= 
\prod_{\alpha_{i} \neq 1}(1 - \alpha_{i}e^{\gamma_{i}}) 
\prod_{\alpha_{i} = 1}(1 - e^{\gamma_{i}}) \\
&= 
\sum_{p \geqslant 0} (-1)^p \op{ch}_{g}(\wedge^p(E_{\neq 0}^{\vee}))
\,\frac{c^{\op{top}}(E_{0})}{\op{Td}(E_{0})}\, ,
\end{split}
\end{equation*}
this together with (\ref{eq_lemme_2}) shows that
\begin{equation}
\label{eq_lemme_3}
\kappa = - c^{\op{top}}(E_{0}) \sum_{j=1}^{r}
\frac{\alpha_{j}e^{\gamma_{j}}}{1 - \alpha_{j}e^{\gamma_{j}}} \, .
\end{equation}
Furthermore, notice that for any function $f(x)$ and 
any $k \in \N$, we have
\[
\frac{d^{k}}{dt^{k}}f(\alpha e^{t}) = \left(
\Big[x \frac{d}{dx}\Big]^{k} f \right)(\alpha e^{t})
\]
and (for $x \neq 1$)
\[
\zeta_{\op{L}}(x,-k) = \Big[x \frac{d}{dx}\Big]^{k} \zeta_{\op{L}}(x,0)
\]
and as well 
\[
\zeta_{\op{L}}(x,0) = \frac{x}{1-x}
\]
we deduce that (for $\alpha \neq 1$)
\begin{equation}
\label{eq_lemme_4}
\frac{\alpha e^{t}}{1 - \alpha e^{t}} = 
\sum_{p \geqslant 0}
\zeta_{\op{L}}(\alpha,-p) \frac{t^{p}}{p !}\, .
\end{equation}
When $\alpha = 1$, we have the the classical expansion
\[
\frac{e^{t}}{1 - e^{t}} = - \frac{1}{t} + 
\sum_{p \geqslant 0}\zeta_{\M{Q}}(-p) \frac{t^{p}}{p!}\, .
\]
Using in (\ref{eq_lemme_3}) the formula just above or the formula (\ref{eq_lemme_4})
according to whether $\alpha_{j}$ is equal to $1$ or not, 
we find that for all $l \geqslant 0$
\[
\kappa^{[l + \op{rk}(E_{0})]} = - c^{\op{top}}(E_{0})
\big[
\sum_{z \in \C{S}_{1}\backslash\{1\}}
\zeta_{\op{L}}(z,-l)\op{ch}^{[l]}(E_{z}^{\vee})
+ \zeta_{\M{Q}}(-l) \op{ch}^{[l]}(E_{0}^{\vee})
\big]
\]
which, noting that $\zeta_{\op{L}}(1,-l) = \zeta_{\M{Q}}(-l)$, concludes the proof.
\endProof

{\bf Caution.} In what follows, in contradiction with classical 
usage and with the Introduction to this article, the notation $L(\chi,s)$ will 
always refer to the {\bf non-primitive} $L$-function associated
with a Dirichlet character $\chi$. We shall write $\chi_{\rm prim}$ for 
the primitive character associated with $\chi$ and accordingly 
write $L(\chi_{\rm prim},s)$ for the associated primitive 
$L$-function.
\smallskip

The following lemma,
proved in \cite[Lemma 5.2, Sec. 5]{LRRIV}, establishes the link between Lerch 
$\zeta$-functions and Dirichlet $L$-functions. It follows from the 
functional equation of Dirichlet $L$-functions when 
the character is primitive. 
 
 \begin{lemma}
 Let $\chi$ be an odd character of $G=\Gal(\mQ(\mn)/\mQ)$.
 The equality
 $$
 \sum_{\sigma\in G}
 \eta(\sigma(\zeta),s)\,\mtr{\chi(\sigma)}
 =
 n^{1-s}{\Gamma(1-s/2)\over\Gamma((s+1)/2)}\pi^{s-1/2}
 L(\mtr{\chi},1-s)
 $$
 holds for all $s\in {\bf C}$.
 \nlabel{SixMonthLem}
 \end{lemma}
 If $\chi$ is a character of $G$, we shall
 write $\tau(\chi):=\sum_{\sigma\in G}\sigma(\zeta)\chi(\sigma)$
 for the Gauss sum associated to $\chi$. Recall that
 if $\chi$ is primitive (i.e. not induced from a subfield
 $\mQ(\mu_m)$ with $m<n$) the following equation holds 
(see \cite[Lemma 4.7, p. 36]{Washington})
\begin{equation}
 \sum_{\sigma\in G}\sigma(\zeta^l)\chi(\sigma)=
 \tau(\chi)\mtr{\chi(l)}
\label{WasEq}
\end{equation}
 where we used the identification $G\simeq(\mZ/n)^\times$ via 
the choice $\zeta=e^{2\pi i/n}$ to give
 meaning to $\chi(l)$.
 
\subsection{Computations}
 
The notations of the sections \ref{introduction} and \ref{preliminaries},
and the conventions of the subsection \ref{invariant}
are still in force. If $N^\times$ is a subgroup of $\mC^\times$ 
and $z,z'\in\mC$, we shall write
 $z\sim_{N^\times} z'$ if $z=\lambda\cdot z'$ with
 $\lambda\in N^\times$. Recall that $f : X \rightarrow \op{Spec}Q_{0}$
is a smooth and projective variety acted upon
by $g$ an automorphism of order $n$ (defined over $Q_{0}$). 
Suppose that $Q_0$ contains $\mQ(\mn)$. 
Let $L$ be a $g$-equivariant ample line bundle over $X$ and 
endow $X(\mC)$ with a K\"ahler metric whose K\"ahler form 
represents the first Chern class of $L$. We will 
denote by $\mtr{\Omega}$ the sheaf of differentials of $f$
equipped with the induced metric. 

We shall now prove Theorems 1 and 2. To do so, 
we first apply the arithmetic fixed point formula
(Theorem \ref{PMainTh}) to the Dolbeault complex $\lambda_{-1}(\mtr{\Omega})$.
 
 \begin{equation*}
 \begin{split}
 f_*(\lambda_{-1}(\mtr{\Omega}))&=
 T_g(\lambda_{-1}(\mtr{\Omega})) \\ 
 &\quad -\int_{X\lmn(\mC)}\Req(TX)\Tdeq(TX)
 \chg(\lambda_{-1}(\Omega))+
 f\umn_*(\lambda_{-1}^{-1}(\mtr{N}^\vee_{X/X\lmn})\lambda_{-1}(\op{\rho}(\mtr{\Omega})))\\
 &=
 T_g(\lambda_{-1}(\mtr{\Omega}))-\int_{X(\mC)_g}\Req(TX)
 \Td(TX_g)\ch(\lambda_{-1}(TX_g^\vee))+
 f_*\umn(\lambda_{-1}(\mtr{\Omega}(f\umn)))\\
 &=
 T_g(\lambda_{-1}(\mtr{\Omega}))-\int_{X(\mC)_g}\Req(TX)c^{\rm top}(TX_g)+
 f\umn_*\bigl(\lambda_{-1}(\mtr{\Omega}(f\umn))\bigr)
 \end{split}
 \end{equation*}
 Applying $\dgm(\cdot)$ to both sides of the last equality, we obtain
 \begin{equation}
 \nlabel{EqMain}
 \begin{split}
 \dgm(f_*(\lambda_{-1}(\mtr{\Omega})))
 &= T_g(\lambda_{-1}(\mtr{\Omega})) \\
 &-{1\over 2}\int_{X(\mC)_g}\Req(TX)c^{\rm top}(TX_g)
 +\dgm(f\umn_*\bigl(\lambda_{-1}(\mtr{\Omega}(f\umn))\bigr))
 \end{split}
 \end{equation}
 We deduce from lemma \ref{torsion} that $T_g(\lambda_{-1}(\mtr{\Omega})) = 0$.
 The following lemma shows that the third term in \refeq{EqMain} likewise vanishes:
 \begin{lemma}
 The equality
 $\dgm(f\umn_*\bigl(\lambda_{-1}(\mtr{\Omega}(f\umn))\bigr))=0$
 holds.
 \end{lemma}
 \beginProof
 Let $d_\mu$ be the relative
 dimension of $X\lmn$ over $\op{Spec}Q_{0}$.
 The expression
 $\dgm(f\umn_{*}\bigl(\lambda_{-1}(\mtr{\Omega}(f\umn))\bigr))$
 can be subdivided into a linear combination of terms of the following kind:
 $$
 \dgm(R^q f\umn_*(\Lambda^{p}(\mtr{\Omega}(f\umn))))+
 \dgm(R^{d_\mu -q}f\umn_*(\Lambda^{d_\mu -p}(\mtr{\Omega}(f\umn))))
 $$
 By Serre duality, the spaces
 $R^q f\umn_*(\Lambda^{p}({\Omega}(f\umn)))$ and
 $R^{d_\mu -q}f\umn_*(\Lambda^{d_\mu -p}({\Omega}(f\umn)))$ are dual
 to each other, and even more, this duality is a duality of
 hermitian vector bundles (for the last statement, see 
\cite{GS2}).
 Hence, from the definition of $\dgm$, it follows that
 $\dgm(R^q f\umn_*(\Lambda^{p}(\mtr{\Omega}(f\umn))))=
 -\dgm(R^{d_\mu -q}f\umn_*(\Lambda^{d_\mu -p}(\mtr{\Omega}(f\umn))))$ which
 permits us to conclude.
 \endProof 
 We shall write $\HDL(\mtr{X})$ for $\HDL(X)$
 equipped with its natural $L^2$-metric. 
 From the preceding discussion and \refeq{EqMain}, we are lead to
 the equality 
 \begin{equation}
 \sum_{k\geqslant 0}(-1)^k\dgm(\HDL^k(\mtr{X}))
 =-{1\over 2}\int_{X(\mC)_g}\Req(TX)c^{\rm top}(TX_g)
 \nlabel{EqMain2}
 \end{equation}
 We shall now show that \refeq{EqMain2} implies
 Theorems 1 and 2. To do so, we will use lemma \ref{SixMonthLem}
 to express derivatives of Lerch $\zeta$-functions
 occuring in $\Req(TX)$ in terms of derivatives of
 Dirichlet $L$-functions, and then lemma \ref{CoHoLe} to give a global (cohomological) expression
 for the right side of \refeq{EqMain2}.
 \smallskip
 
 {\bf Proof of Theorem 1.}
 We compute
 $$
 \sum_{k\geqslant 0}(-1)^k\,\dgm(\HDL^k(\mtr{X})) =
 \sum_{k\geqslant 0}(-1)^k\sum_{l\geqslant 0}\zeta^l\,\dgm(\HDL^k(\mtr{X})_l)
 $$
 and
 
 \begin{equation*}
 \begin{split}
 \sum_{\sigma\in G}\sum_{k\geqslant 0}(-&1)^k\sum_{l\geqslant 0}
 \sigma(\zeta)^l\,\dgm(\HDL^k(\mtr{X})_l)\chi(\sigma)\\
 &=
 \tau({\chi})\sum_{k\geqslant 0}
 (-1)^k\sum_{l\geqslant 0}\mtr{\chi(l)}\op{\dgm}(\HDL^k(\mtr{X})_l)\\
 &=
 -{\tau({\chi})}\int_{X(\mC)_g}
 L\pr(\mtr{\chi},0)\sum_l\mtr{\chi}(l)\op{rk}(TX_l)c^{\rm top}(TX_g)\\
 &=
 -{\tau({\chi})}
 \int_{X(\mC)_g}{L\pr(\mtr{\chi},0)\over L(\mtr{\chi},0)\tau({\chi})}
 \sum_{\sigma\in G}\chi(\sigma)\Td(TX_g)
 {\sum_{p\geqslant 0}(-1)^p p\cdot\ch_{\sigma(\zeta)}(\Lambda^p(\Omega))
 \over
 \sum_{p\geqslant 0}(-1)^p\ch_{\sigma(\zeta)}(\Lambda^p(N^\vee))}\\
 &=
 -{\tau(\chi)}{L\pr(\mtr{\chi},0)\over L(\mtr{\chi},0)}
 \sum_{\sigma\in G}\sum_{p,q}(-1)^{p+q}p\cdot\rk(H^{p,q}(X(\mC))_\sigma)
 \chi(\sigma)
 \end{split}
 \end{equation*} 
 which is the result. For the first equality in the last string of equalities, we have 
used \refeq{WasEq}; for the second one, we have used \refeq{EqMain2} and applied 
\refeq{WasEq} to the definition of Lerch's partial $\zeta$-functions; for the third one, 
we used lemma \ref{CoHoLe};  
for the last equality, we have 
 applied the holomorphic Lefschetz trace formula \cite[3.4, p. 422]{GH} 
to the virtual vector bundle 
$1-\Omega+2\cdot\Lambda^2(\Omega)-\dots+(-1)^{\dim(X)}\dim(X)\cdot
\Lambda^{\dim(X)}(\Omega)$. 
To conclude the proof, we still need the following lemma.  
Write $\mot:=\C{X}(X,g)$.
\begin{lemma}
Let $\sigma\in\Hom(\mQ(\mn),\mC)$.
If $\omega_{\sigma}\in \det_\mC H(\mot,\mC)_\sigma$ is non-zero and defined over $Q_0$ 
for the de Rham $Q_0$-structure of $H(\mot,\mC)$, then 
the equality $|P_\sigma(\mot)|^2 \sim_{|Q_0^\times|} |\omega_{\sigma}|_{L^2}^2$ holds.
\label{PerLem}
 \end{lemma}
 \beginProof
% Before beginning with the proof, notice the following. If we equip 
% $\det_\mC(H(\mot)_\sigma)$ and $\otimes_{k=1}^{\rk_\mQ
% (H(\mot))}H(\mot)_\mC$ with their natural metrics, the 
%embedding $\det_\mC(H(\mot)_\sigma)\hookrightarrow \otimes_{k=1}^{\rk_\mQ
%(H(\mot))}H(\mot)_\mC$ given by ... is an isometry onto its image.  
Let $z$ be a complex number such that 
$z\cdot \omega_{\sigma}$ is defined over the singular $Q_0$-structure of 
$\det_\mC H(\mot,\mC)_\sigma$. Since $z\cdot \omega_{\sigma}$ is of pure 
Hodge type, the construction of the Hodge 
metric and lemma \ref{HDRLem} shows that the
number $|z\cdot \omega_{\sigma}|_{L^2}$ lies in $|Q_0^{\times}|$ and from 
this we can conclude.
 \endProof

 {\bf Proof of Theorem 2.}
 By class-field theory and proposition \ref{BasInv}, we are reduced to
 the case $Q=\mQ(\mn)$. We may thus suppose that $X$ is an abelian 
 variety $A$ with (not necessarily
 maximal) complex multiplication by ${\cal O}_{\mQ(\mn)}$. 
 Let $\mtr{V}:=\HDL^1(\mtr{A})$.
 By the usual Lefschetz trace formula, the number of fixed points 
on $A(\mC)$ of the automorphism corresponding to $\zeta\in 
{\cal O}_{\mQ(\mn)}$ 
  is $\prod_{l\in\Zn}(1-\zeta^l)^{\rk(V_l)}$ and therefore 
we get, using \refeq{EqMain2} and lemma \ref{2Ten}
 $$
 -\sum_l{\zeta^l\over 1-\zeta^l}\dgm(\mtr{V}_l)=
 -{1\over 2}R_g(TA|_0)=
 -{1\over 2}\sum_l 2i(\partial/\partial s)\eta(\arg(\zeta^l),0)
 \rk((TA|_0)|_l).
 $$
We shall write $TA_l$ for $(TA|_0)_l$. 
 Now notice that $\NIm{\zeta^l\over 1-\zeta^l}=
 \eta(\arg(\zeta^l),0)$ and thus
 $$
 2\sum_l\eta(\zeta^l,0)\dgm(\mtr{\Omega}_l)=
 \sum_l(\partial/\partial s)\eta(\arg(\zeta^l),0)\rk(TA_l).
 $$
Next we take the Fourier transform of both sides of the last equality   
for the action of $G=\Gal(\mQ(\mn)/\mQ)$. We get:
 $$
 -2\sum_l[\sum_{\sigma\in G}\eta(\arg(\sigma(\zeta)^l),0)\chi(\sigma)]
 \dgm(\mtr{\Omega}_l)=
 \sum_l[\sum_{\sigma\in G}(\partial/\partial
 s)\eta(\arg(\sigma(\zeta)^l),0)\chi(\sigma)]
 \rk({\Omega}_l)
 $$
 and changing variables 
\begin{multline*}
 -2[\sum_{\sigma\in G}\eta(\arg(\sigma(\zeta)),0)\chi(\sigma)]
 \sum_l\mtr{\chi}(l)\dgm(\mtr{\Omega}_l)   \\
 =
 [\sum_{\sigma\in G}(\partial/\partial s)\eta(\arg(\sigma(\zeta)),0){\chi}(\sigma)]
 \sum_l\mtr{\chi}(l)\rk(\Omega_l).
\end{multline*}
 We now calculate, using lemma \ref{SixMonthLem}
 \begin{equation*}
 \begin{split}
 \frac{[\sum_{\sigma\in G}(\partial/\partial 
s)\eta(\arg(\sigma(\zeta)),0){\chi}(\sigma)]}
 {[\sum_{\sigma\in G}\eta(\arg(\sigma(\zeta)),0)\chi(\sigma)]}&\\
 =
 -\log(n)-{1\over 2}({\Gamma\pr(1)\over \Gamma(1)}&+
 {\Gamma\pr(1/2)\over \Gamma(1/2)})+
 \log(\pi)-{L\pr({\chi},1)\over L({\chi},1)}\\
 =
 -\log(n)-{1\over 2}({\Gamma\pr(1)\over \Gamma(1)}&+
 {\Gamma\pr(1/2)\over \Gamma(1/2)})+
 \log(\pi) \\
 & - {L\pr({\chi_\prim},1)\over L({\chi_\prim},1)}-
 \sum_{p|n}{{\chi_\prim}(p)\over
 p-{\chi_\prim}(p)}\log(p)\\
 = -\log(n)-{1\over 2}({\Gamma\pr(1)\over \Gamma(1)}&+
 {\Gamma\pr(1/2)\over \Gamma(1/2)})+
 \log(\pi)-\log({2\pi\over f_\chi})\\
 &+
 {\Gamma\pr(1)\over\Gamma(1)}+
 {L\pr(\mtr{\chi}_\prim,0)\over L(\mtr{\chi}_\prim,0)}-
 \sum_{p|n}{{\chi_\prim}(p)\over
 p-{\chi_\prim}(p)}\log(p)
 \end{split}
 \end{equation*}
 \begin{equation*}
 =
 \log({f_\chi\over n})+{L\pr(\mtr{\chi}_\prim,0)\over
 L(\mtr{\chi}_\prim,0)}-\sum_{p|n}{{\chi_\prim}(p)\over
 p-{\chi_\prim}(p)}\log(p).\quad\quad\quad\quad \;\;
 \end{equation*}
Together with lemma \ref{PerLem}, this allows us to conclude. 
We used the functional equation of primitive 
Dirichlet $L$-functions for the third equality. 
 \endProof

 \section{The period conjecture of Gross-Deligne}

 In this section, we shall indicate the consequences of
Theorem 1 and Theorem 2 for the period conjecture of Gross-Deligne 
\cite[Sec. 4, p. 205]{Gross}.
 We first recall the latter conjecture.
 Let $Q$ be a finite abelian extension of $\mQ$ and
 let $H$ be a rational and homogenous Hodge structure of dimension $[Q:\mQ]$.
 Suppose that there is a morphism of rings
 $\iota:Q\hookrightarrow\End(H)$ (in other
 words, $H$ has maximal complex multiplication
 by $Q$). Suppose also
 that $H$ is embedded in the singular cohomology
 $H(X,\mQ)$ of a variety $X$ defined over $\Qb$.
 We let $f_{Q}$ be the conductor of $Q$
 and we choose an embedding of $Q$ in $\mQ(\mu_{f_{Q}})$ 
(this is possible by class-field theory). Choose an embedding
 $\varphi:Q\hookrightarrow\mC$ and an isomorphism
 $\Gal(\mQ(\mu_{f_{Q}})/\mQ)\simeq(\mZ/f_{Q})^{\times}$.
 There is a natural map $\Gal(\mQ(\mu_{f_{Q}})/\mQ)\ra
 \Hom(Q,\mC)$ given for each $\sigma\in\Gal(\mQ(\mu_{f_{Q}})/\mQ)$ by $\varphi\circ\sigma|_Q$,
 and we thus obtain a map $(\mZ/f_{Q})^{\times}\ra \Hom(Q,\mC)$.
 For each $u\in(\mZ/f_{Q})^{\times}$ let
 $\omega_{u}^{H}\in H\otimes_\mQ\mC$ be a 
 non-vanishing element affording the embedding corresponding to $u$ 
 and defined over $\Qb$ for the de Rham 
 $\Qb$-structure of $H(X,\mC)$
 (such an element is well-defined up to multiplication by a
 non-zero algebraic number). We attach to
 $\omega_u^H$ a period $\Per(\omega_u^H):= v(\omega_u^H)$ where $v\in H^\vee$ is
 any (non-zero) element of the dual of the $\mQ$-vector space $H$.
 The number $\Per(\omega_u^H)$
 is independent of the choices of $v$ and $\omega_u^H$,
 up to multiplication by a
 non-zero algebraic number,
 and only depends on $u$ and $H$. 
 In the notation of the beginning of 
 the Introduction, with $Q_0=\Qb$, we have 
 $\Per(\omega_u^H)=1/P_u(H)$ (where we identify $u$ with the corresponding 
 embedding of $Q$ after the equality sign).
 Let $(p(u),q(u))$ be the Hodge type of
 $\omega_{u}^{H}$. By \cite[Lemme 6.12]{Deligne2},
 there exists a (non-unique) function
 $\epsilon^H:\mZ/f_{Q}\ra\mQ$ which satisfies the equation
 \begin{equation*}
 p(u)=\sum_{a\in\mZ/f_{Q}}\epsilon^H(a)[u\cdot a/f_{Q}]
% \nlabel{SysDel}
 \end{equation*}
 for all $u\in(\mZ/f_{Q})^{\times}$. Here $[\cdot]$ takes
 the fractional part. The following conjecture is
 formulated by Gross in \cite[p. 205]{Gross}; he 
indicates that the
 precise form of it was suggested to him by Deligne.
 This conjecture is related by Deligne to his conjecture
 on motives of rank 1 in \cite[8.9, p. 338]{Deligne}. 
 \medskip

 {\bf Period conjecture.}
 {\it Let $u\in(\mZ/f_{Q})^\times$. The relations
 
~
\begin{equation*}
\tag{i}
 |\Per(\omega_u^H)|\sim_{\Qb^\times}\prod_{a\in \mZ/f_{Q}}
 \Gamma(1-{a\over f_{Q}})^{\epsilon^H(a/u)}
\end{equation*}
~
\begin{equation*}
\tag{ii}
\Per(\omega_u^H)\sim_{\Qb^\times}|\Per(\omega_u^H)|
\end{equation*}
 hold.
 \nlabel{conjP1}}
 \medskip\\
The relations (i) and (ii) can of course be condensed 
 in the single relation $\Per(\omega_u^H)
 \sim_{\Qb^\times}\prod_{a\in \mZ/f_{Q}}
 \Gamma(1-{a\over f_{Q}})^{\epsilon^H(a/u)}$.
   
 One can show that the period conjecture is independent of the choice
 of the function $\epsilon^H$ (see the appendix 
by Koblitz and Ogus to \cite{Deligne}).\\
We let $\Gamma_{f_{Q}}$ be the algebraic extension of 
$\mQ$ generated by the numbers 
$$\prod_{a\in \mZ/f_{Q}}\Gamma(1-{a\over f_{Q}})^{\gamma(a/u)},$$ 
where 
$\gamma$ runs on all the functions
 $\mZ/f_{Q}\ra\mQ$ satisfying the equation
 \begin{equation*}
 \sum_{a\in\mZ/f_{Q}}\gamma(a)[u\cdot a/f_{Q}]=0
 \end{equation*}
 for all $u\in(\mZ/f_{Q})^{\times}$.

 Suppose now that $Q=\mQ(\mu_{p})$ for 
 a prime number $p$ and let $Q_0\subseteq\Qb$ be a field 
of definition of $X$ containing
$Q$ and $\Gamma_{p}$. Suppose that $H$ inherits 
the de Rham $Q_0$-structure of $H(X,\mC)$ and that 
its $Q$-vector space structure is compatible 
with that $Q_0$-structure. Suppose furthermore that all 
 the $\omega_u^H$ are defined over $Q_0$ 
 (for the de Rham $Q_0$-structure) and that
 $\Per(\omega_u^H)\sim_{Q_0^\times}{1\over\Per(\omega_{-u}^H)}$ for all $u$.\\
If $\chi$ is an Artin character, we define $\mQ(\chi)$ as the field 
generated over $\mQ$ by the values of $\chi$. 
We denote by $E$ the compositum in $\mC$ 
of the field $\mQ(\mu_{p})$ and of the fields $\mQ(\chi)$ for 
all the odd Artin characters $\chi$ of $\mQ(\mu_p)$, and we let
$E^+:=E\cap\mR$. 
\begin{lemma}
 If the conjecture ${\rm B}(H,E,\chi)$ holds 
for all the odd Artin characters $\chi$ of 
$Q=\mQ(\mu_p)$, then the identity 
 $$
 \log|\Per(\omega_u^H)|=\log|\prod_{a\in \mZ/p}
 \Gamma(1-{a\over p})^{\epsilon^H(a/u)}|+\sum_i b_{u,i}\log|a_{u,i}|
 $$
holds for all $u$. Here $b_{u,i}\in E^+$, $a_{u,i}\in Q_0$ and  
$i$ runs over a finite
 set of indices.
\label{LinkLem}
\end{lemma}
 \beginProof
We can plainly assume that $\epsilon^H(0)=0$. 
 We write
 \begin{eqnarray*}
 \lefteqn{
 \sum_{u}(\log|\Per(\omega_u^{H})|-\sum_i b_{u,i}\log|a_{u,i}|)\chi(u)=
 \sum_{u}\chi(u)\sum_{a}\log(\Gamma(1-a/p))\epsilon(a/u)}
 \\&=&
 \sum_a\log(\Gamma(1-a/p))\sum_u\chi(u)\epsilon(a/u)=
 \sum_a\log(\Gamma(1-a/p)\chi(a)\sum_u\chi(u)\epsilon(1/u).
 \end{eqnarray*}
Now by the definition of $\epsilon$
 $$
 \sum_u\chi(u)p(u)=
 \sum_a[a/p]\sum_u\chi(u)\epsilon(a/u)=
 (\sum_a\chi(a)[a/p])(\sum_u\chi(u)\epsilon(1/u))
 $$
 and thus
 $$
 \sum_u(\log|\Per(\omega_u^{H})|-\sum_i b_{u,i}\log|a_{u,i}|)\chi(u)=
 (\sum_a\log(\Gamma(1-a/p))\chi(a)){\sum_u\chi(u)p(u)\over
 \sum_a\chi(a)[a/p]}
 $$
 By Hurwitz's formula, this equals
\begin{equation*}
\begin{split}
({L\pr(\chi,0)\over L(\chi,0)}&+\log(p))\sum_u\chi(u)p(u)\\
&=
{L\pr(\chi_\prim,0)\over L(\chi_\prim,0)}\sum_u\chi(u)p(u) +
(\log(p)+{\chi_\prim(p)\over 
1-\chi_\prim(p)}\log(p))\sum_u\chi(u)p(u)
\end{split}
\end{equation*}
which is the result.
 \endProof

Suppose moreover that $X$ is acted upon 
by an automorphism $g$ (defined over $Q_{0}$) 
of order $p$, and let
$d_k:=\dim_{Q}(H^k(\C{X}(g),\mQ))$ for any $k\geqslant 0$. 
Arising from the K\"unneth isomorphism 
$H(\times_{j=1}^{d_k} X,\mQ)\simeq\otimes_{j=1}^{d_k}H(X,\mQ)$ and 
lemma \ref{EmbHodge},
 there is a natural embedding of Hodge structures
 $\det_{Q}(H^k(\C{X}(g),\mQ))\hookrightarrow
 H(\times_{j=1}^{d_k} X,\mQ)$ which respects the de Rham 
$Q_0$-structures. 

The next two results
follow from Theorem 1 and lemma \ref{LinkLem}.
\begin{cor}
 Suppose that there is at most one $k$ such that $0\leqslant k\leqslant \dim(X)$ 
and $d_k\not = 0$. The 
 natural embedding of Hodge structures
 $\det_{Q}(H^k(\C{X}(g),\mQ))\hookrightarrow
 H(\times_{j=1}^{d_k} X,\mQ)$ satisfies the hypothesis of 
the period conjecture and the identity 
 $$
 \log|\Per(\omega_u^{\det_Q(H^k(\C{X}(g),\mQ))})|=
\log|\prod_{a\in \mZ/p}
 \Gamma(1-{a\over p})^{\epsilon^{\det_Q(H^k(\C{X}(g),\mQ))}(a/u)}|+\sum_i 
b_{u,i}\log|a_{u,i}|
 $$
holds for all $u$, 
where $b_{u,i}\in E^+$, $a_{u,i}\in Q_0$ and $i$ runs over a finite
 set of indices.
\end{cor}
Notice that in particular the last corollary applies if $X$ is a 
hypersurface of a projective space which carries an automorphism 
extending $g$. 
\begin{cor}
 If $X$ is a surface, the natural embedding of Hodge structures 
 $\det_{Q}(H^2(\C{X}(g),\mQ))\hookrightarrow
 H(\times_{j=1}^{d_2} X,\mQ)$ satisfies the hypothesis of 
the period conjecture and the identity 
 $$
 \log|\Per(\omega_u^{\det_Q(H^2(\C{X}(g),\mQ))})|=
\log|\prod_{a\in \mZ/p}
 \Gamma(1-{a\over p})^{\epsilon^{\det_Q(H^2(\C{X}(g),\mQ))}(a/u)}|+\sum_i 
b_{u,i}\log|a_{u,i}|
 $$
holds  for all $u$, 
where $b_{u,i}\in E^+$, $a_{u,i}\in Q_0$ and $i$ runs over a finite
 set of indices.
\end{cor}
\beginProof
 The conjecture ${\rm B}(H^1(\C{X}(g),\mQ),E,\chi)$ is verified for 
all the odd Artin characters $\chi$ of $\mQ(\mu_p)$ because of the existence of
 the Picard variety and Theorem 2. The lemma \ref{LinkLem} 
 now implies that the conjecture is verified
 for $H^2(\C{X}(g),\mQ)$.
\endProof
Notice that when $p=3$ 
 in either of the last corollaries, 
the item (i) of the period 
conjecture holds for the described embedding of Hodge structures 
(for all $u\in(\mZ/p)^\times$), 
since $E^+=\mQ$ in that case.\\
{\bf Remark.} The (ii) of the period conjecture seems to be out of the 
reach of the techniques developped in the present paper.

 \end{document}